\documentclass[12pt]{article}
\usepackage[top=1in, bottom=1in, left=0.5in, right=0.5in]{geometry}
\usepackage{amssymb}
\usepackage{amsmath}
\usepackage{amsthm}
\usepackage{mathrsfs}
\usepackage[colorlinks=true, allcolors=black]{hyperref}
\usepackage{graphicx}
\numberwithin{equation}{section}
\title{\textbf{{\Large Duality Between Prime Factors and The
Prime Number Theorem For Arithmetic Progressions - Higher Order Dualities}}}

\author{\textit{Krishnaswami Alladi} and \textit{Sroyon Sengupta}}
\date{}

\begin{document}

\maketitle

\textbf{Abstract:} \textit{{\small In 1977, the first author observed a duality between the largest and smallest prime factors of integers, and established as a consequence some new results on the Möbius function $\mu(n)$ using the Prime Number Theorem for Arithmetic Progressions. In that 1977 paper, higher order dualities were observed involving the $k$-th largest and $k$-th smallest prime factors, facilitated by the M\"obius function and $\omega(n)^{k-1}$, where $\omega(n)$ is the number of distinct prime factors on $n$. In 2024, the first author and Jason Johnson proved new results involving $\mu(n)$ and $\omega(n)$, by exploiting the second order duality identity of Alladi (1977). We establish here extensions to all higher orders $k$, the results of Alladi (1977) and of Alladi-Johnson (2024), by utilizing the $k$-th order duality in Alladi's 1977 paper. First, we show that for each $k\geq 2$,}}
\begin{eqnarray}
        \sum_{n=2}^{\infty} \frac{\mu(n)\omega(n)^{k}}{n} =0, \nonumber
\end{eqnarray}
    \textit{{\small where $\mu(n)$ is the M\"obius Function and $\omega(n)$ counts the number of distinct prime factors of $n$.   Further, using the General Duality Identity and the Prime Number Theorem of Arithmetic  Progressions, we prove that for integers $j,\ell$ satisfying $1 \leq j \leq \ell$ and $(j,\ell)=1$}} 
    \begin{eqnarray}
        \sum_{\substack{n=2 \\ p_1(n) \equiv j\;(mod\;\ell)}}^{\infty} \frac{\mu(n)\omega(n)^{k-1}}{n}=0, \nonumber
\end{eqnarray}
\textit{{\small for every $k \geq 3$; this result for $k=1$ is due to Alladi (1977) and for $k=2$ due to Alladi-Johnson (2024). We also recast this result in the following manner as a density-type theorem: for integers $j,\ell$ satisfying $1 \leq j \leq \ell$ and $(j,\ell)=1$}}
    \begin{eqnarray}
      (-1)^k\sum_{\substack{n=2 \\ p_1(n) \equiv j\;(mod\;\ell)}}^{\infty} \frac{\mu(n){\omega(n)-1 \choose k-1}}{n}=\frac{1}{\varphi(\ell)}, \nonumber
      \end{eqnarray}
    \textit{{\small for every $k \geq 3$. All results are established here in quantitative form.}} 

\section{Introduction and Background}
Two famous results by Edmund Landau \cite{LaT99} are that 
\begin{eqnarray}
    M(x):=\sum_{1 \leq n \leq x}\mu(n) = o(x),\quad \text{as $x \rightarrow \infty$} 
\end{eqnarray}
and
\begin{eqnarray}
    \sum_{n=1}^{\infty} \frac{\mu(n)}{n} = 0 
\end{eqnarray}
are elementarily equivalent to the Prime Number Theorem (PNT), where $\mu(n)$ is the M\"obius Function. There are also results equivalent to the Prime Number Theorem in Arithmetic Progressions (PNTAP), with $\mu(n)$ is replaced by $\mu(n)\chi(n)$, where $\chi(n)$ is the Dirichlet character modulo $k$, when the arithmetic progression considered in context has a common difference $k$. \\ \\
In 1977 \cite{KA1977}, the first author noticed the following interesting Duality identities involving
the M\"obius function that connect the smallest and largest prime factors: 
\begin{eqnarray}
    \sum_{2 \leq d|n}\mu(d)f(p_1(d)) = -f(P_1(n)) \;\;\;\text{and}\;\;\; \sum_{2 \leq d|n}\mu(d)f(P_1(d)) = -f(p_1(n)),
\end{eqnarray}
where, for integers $n>1$, $p_1(n)$ and $P_1(n)$ are the smallest and largest prime factors of $n$ respectively, and $f$ is ANY function defined on primes. Using the first identity in (1.3) along with properties of the M\"obius function, it was proved in \cite{KA1977} that if $f$ is a bounded function on the primes satisfying
\begin{eqnarray}
    \lim_{x \to \infty} \frac{1}{x} \sum_{1 \leq n \leq x} f(P_1(n)) = c,
\end{eqnarray}
then
\begin{eqnarray}
    \sum_{n =2}^{\infty}\frac{\mu(n)f(p_1(n))}{n} = -c, 
\end{eqnarray}
and vice-versa. This can be realized as a generalization of Landau's result as follows. We consider the function $f$ such that $f(p)=1$ for all primes $p$. Then $c$ in (1.4) is 1. Further, rewrite (1.2) as
\begin{eqnarray}
    \sum_{n=2}^{\infty} \frac{\mu(n)}{n}=-1.
\end{eqnarray}
This is exactly what we have in (1.5) for the chosen $f$. \\ \\
It was shown in \cite{KA1977} that PNTAP implies that the sequence $\{P_1(n)\}$ of the largest prime factors of $n$ is uniformly distributed in the reduced residue classes modulo a positive integer $\ell$. In other words, if $f$ is chosen to be the characteristic function of the primes in a given arithmetic progression $j\;(mod\;\ell)$, then for such an $f$, equation (1.4) is satisfied with $c = \frac{1}{\varphi(\ell)}$, and therefore, by (1.5), we have
\begin{eqnarray}
    \sum_{\substack{n =2 \\ p_1(n) \equiv j\;(mod\; \ell)}}^{\infty} \frac{\mu(n)}{n} = -\frac{1}{\varphi(\ell)},
\end{eqnarray}
for ALL positive integers $\ell$ and $1 \leq j < \ell$ satisfying $(j,\ell)=1$. This is quite an intriguing result since the LHS of (1.7) is a subseries of what we have in (1.6), and so each of those $\varphi(\ell)$ many slices (due to $\varphi(\ell)$ many possible $j$'s) have the exact same value, and these add up to give $-1$. In the past decade, (1.7) has been extended and generalized in the setting of algebraic number theory by several authors \cite{Da17}, \cite{SW19}, \cite{Wa20}, and \cite{KMS}.  \\ \\
In \cite{KA1977}, the following general duality identities were also noted: For a positive integer $k$, let $P_k(n)$ and $p_k(n)$ denote the $k^{th}$ largest and smallest prime factors of $n$, respectively, if $n$ has at least $k$ distinct prime factors. Let $\omega(n)$ denote the number of distinct prime factors of $n$. Then:
\begin{eqnarray}
    \sum_{1<d|n}^* \mu(d)f(P_k(d)) = (-1)^k {\omega(n)-1 \choose k-1}f(p_1(n)),
\end{eqnarray}
and
\begin{eqnarray}
    \sum_{1<d|n}^*\mu(d)f(p_k(d)) = (-1)^k {\omega(n)-1\choose k-1}f(P_1(n)),
\end{eqnarray}
where $*$ over the summation means that if $n$ has fewer than $k$ distinct prime factors, then the sum is zero.  Now, by M\"obius Inversion, it follows from (1.8) and (1.9) that
\begin{eqnarray}
    \sum_{1<d|n} \mu(d){\omega(d)-1 \choose k-1}f(p_1(d)) = (-1)^kf(P_k(n)),
\end{eqnarray}
and
\begin{eqnarray}
    \sum_{1<d|n}\mu(d){\omega(d)-1 \choose k-1}f(P_1(d)) = (-1)^kf(p_k(n)).
\end{eqnarray}
Note here that we have omitted the $*$ by adopting the convention that $f(P_k(n))=0=f(p_k(n))$, whenever $\omega(n)<k$, i.e. if $n$ has fewer than $k$ distinct prime factors (in that case, the condition represented by $*$ earlier stays consistent). \\ \\
Motivated by the Second Order Duality, namely the case $k=2$ in (1.8) and (1.10), the first author along with Jason Johnson proved the following in \cite{AJ24}: for integers $j,\ell$ satisfying $1\leq j \leq \ell$ and $(j,\ell)=1$
\begin{eqnarray}
    \sum_{\substack{n=2 \\ p_1(n) \equiv j \;(mod\;\ell)}}^{\infty} \frac{\mu(n)\omega(n)}{n} = 0.
\end{eqnarray}
In the course of proving (1.12), several essential auxiliary results were proved, a crucial one being the uniform distribution of the second largest prime factor $P_2(n)$ in the reduced residue classes modulo a positive integer $\ell$. It is also proven in \cite{AJ24} that
\begin{eqnarray}
    \sum_{n=2}^{\infty} \frac{\mu(n)\omega(n)}{n}=0 .
\end{eqnarray}
Therefore, quite similar to that of the results in \cite{KA1977}, we see that slicing up the sum in (1.13) into $\varphi(\ell)$ many partial sums, we still get that each of them are respectively 0, adding up to 0. \\ \\
A general theorem later noted in \cite{AJ24} is that if $f$ is a bounded function on the primes such that 
\begin{eqnarray}
    \sum_{2 \leq n \leq x} f(P_1(n)) \sim \kappa x
\end{eqnarray}
and
\begin{eqnarray}
    \sum_{2 \leq n \leq x} f(P_2(n)) \sim \kappa x
\end{eqnarray}
for some constant $\kappa$, then
\begin{eqnarray}
    \sum_{n=2}^{\infty} \frac{\mu(n)\omega(n)f(p_1(n))}{n}=0.
\end{eqnarray}
Of course, considering $f$ to be the characteristic function of the primes in a given arithmetic progression $j \;(mod\;\ell)$, then for such an $f$, $\kappa = \frac{1}{\varphi(\ell)}$, due to the uniform distribution of the first and the second largest prime factors of $n$ in reduced residue classes, and thence, we get the original result (1.12). \\ \\
In this paper, we establish extensions of all results proved in \cite{AJ24} for all integers $k \geq 3$. We first start by proving that for any arbitrary positive integer $k \geq 2$
\begin{eqnarray}
    \sum_{n=2}^{\infty} \frac{\mu(n)\omega(n)^{k}}{n} = 0.
\end{eqnarray}
We actually establish this in strong quantitative form by providing a series expansion for the sum in (1.17) truncated at $x$ (see \textit{Theorem 2.1-k}) in \S 2. We then move on to prove that for any arbitrary positive integer $k$, the sequence of the $k^{th}$ largest prime factors $\{P_k(n)\}$ of $n$ is uniformly distributed in reduced residue classes (see \textit{Theorem 5.3} in \S5). Using this, and the duality identity (1.10), we establish our main result (\textit{Theorem 6.5} \S6): \textit{Given any fixed but arbitrary positive integer $k \geq 3$, and for integers $j,\ell$ satisfying $1\leq j \leq \ell$ and $(j,\ell)=1$, we have that}
\begin{eqnarray}
    \sum_{\substack{n=2 \\ p_1(n) \equiv j \;(mod\;\ell)}}^{\infty} \frac{\mu(n)\omega(n)^{k-1}}{n} = 0.
\end{eqnarray}
We also give a density type result that states as follows:  \textit{Given any fixed but arbitrary positive integer $k \geq 3$, and for integers $j,\ell$ satisfying $1\leq j \leq \ell$ and $(j,\ell)=1$, we have that}
\begin{eqnarray}
     (-1)^k\sum^{\infty}_{\substack{n=2 \\ p_1(n)\equiv j\;(mod\;\ell)}} \frac{\mu(n){\omega(n)-1 \choose k-1}}{n} = \frac{1}{\varphi(\ell)}.
\end{eqnarray}
All our results are proved in quantitative form. \\ \\
\textit{\underline{Remark:}} We will first prove (1.18) and (1.19) for $k=3$. The proofs for every $k \geq 4$ will follow quite similarly, but the asymptotics get more involved as $k$ increases. So we provide the full details for $k=3$, and for $k\geq 4$, we mention only the main steps of the proof.  
%So we focus on proving results in detail for the case $k=3$, which is where certain features of the asymptotics are different from the case $k=2$ investigated in [A-J 2024]. The situation regarding $k\geq 3$, although increasing in complexity as $k$ increases, is similar to $k=3$,
%and so we only sketch how the methods used to study the case $k=3$ yield corresponding results for each $k\ge 4$.

\subsection{Notations and Conventions}
We denote by $c_1,c_2,\cdots  $ absolute positive constants, whose values will not concern us. The $\ll$ and $O$ notations are equivalent and will be used interchangeably, as is convenient. We adopt the convention that 
\begin{eqnarray}
    f(x) \ll g(x) \implies |f(x)| \leq Kg(x),
\end{eqnarray}
with $x$ ranging in some domain depending on the context, and $K$ a positive constant. Implicit constants are absolute unless otherwise indicated with a subscript. Although our results can be established with uniformity by allowing the modulus $\ell$ to grow slowly as a function $x$, we prove the results by fixing an arbitrary $\ell$. The alphabet $n$, whether used as an argument of a function, or in a summation, will always be a positive integer, unless specified otherwise. Further, $p,q,r$ and $p_i$'s for positive integer subscripts $i$, will always represent prime numbers, whether used in an argument of a function or a summation. The parameter $T=T(x)$, a function that will come up in most of the proofs, will be chosen optimally to get suitable and required bounds in various estimates, and may not be the same in all contexts. We shall use the standard notation $[x]$ for the integral part of a real number $x$, and $\{x\}$, where indicated, will denote the fractional part of $x$. Thus, $\{x\}=x-[x]$. Finally, complex numbers will be denoted either by $z$, or by $s=\sigma+it$, the latter while dealing with Dirichlet series. Further notation will be introduced in the sequel as needed. \\ \\
\textbf{Note:} Whenever we employ the duality identity of order $k$, we will use $\omega(n)^{k-1}$ in view of (1.10). Otherwise, we will simply use $\omega(n)^k$. 

\section{Sums of the M\"obius function and $\omega(n)^k$}
The aim of this section is to prove that for all positive integers $k \geq 2$,
\begin{eqnarray}
    M_{\omega^k}(x):=\sum_{n\leq x}\mu(n)\omega(n)^k = o(x),
\end{eqnarray}
and
\begin{eqnarray}
    m_{\omega^k}(x):= \sum_{n\leq x}\frac{\mu(n)\omega(n)^k}{n} = o(1), \quad \text{as $x \to \infty$}.
\end{eqnarray}
We obtain strong quantitative versions of (2.1) and (2.2) in the form of an asymptotic estimate given by a series in decreasing powers of $\log x$ for $M_{\omega^k}(x)$, and a similar asymptotic estimate for $m_{\omega^k}(x)$ given by a series in decreasing powers of $\log x$ as well (see \textit{Theorem 2.1} and \textit{2.1-k} in this section). \\ \\
For the case $k=1$, Alladi-Johnson established (2.1) and (2.2) by obtaining bounds and communicated this to Tenenbaum. In response to this, Tenenbaum in a letter to Alladi stated that using the Selberg-Delange analytic method, the following sharper estimates hold: \textit{If $\nu$ is an arbitrary but fixed positive integer, then there exist constants $\beta_j$ such that}
\begin{eqnarray}
    M_{\omega}(x) = \frac{x}{\log^2x}+\frac{\beta_3x}{\log^3x}+\frac{\beta_4x}{\log^4x}+\cdots  +\frac{\beta_{\nu}x}{\log^{\nu}x}+O_{\nu}\left(\frac{x\log\log^{2\nu+2}x}{\log^{\nu+1}x}\right),
\end{eqnarray}
\textit{and there exist constants $\gamma_j$ such that}
\begin{eqnarray}
    m_{\omega}(x) = -\frac{1}{\log x}+\frac{\gamma_2}{\log^2x}+\frac{\gamma_3}{\log^3x}+\cdots  +\frac{\gamma_{\nu}x}{\log^{\nu}x}+O\left(\frac{\log\log^{2\nu+4}x}{\log^{\nu+1}x}\right).
\end{eqnarray}
Tenenbaum's error terms were in a slightly different form, but (2.3) and (2.4) are essentially equivalent to what he stated. In \cite{AJ24}, Alladi-Johnson proved (2.3) and (2.4) by an elementary method, but using the Prime Number Theorem in the following strong form \begin{align}
    \pi(x) &= \ell i(x)+O(xe^{-c\sqrt{\log x}})\nonumber \\
            &= \frac{x}{\log x}+\frac{x}{\log^2x}+\frac{2!x}{\log^3x}+\cdots  +\frac{\log^{\nu}x}{(\nu-1)! x}+O\left(\nu!\frac{x}{\log^{\nu +1}x}\right),
\end{align}
given by repeated integration-by-parts of $\ell i(x)$. Here, we shall follow the elementary method of Alladi-Johnson to establish \textit{Theorems 2.1} and \textit{2.1-k} below for fixed $k \geq 2$, but we shall rely on strong asymptotic estimates for
\begin{eqnarray}
    \pi_k(x):=\sum_{\substack{n \leq x \\ \omega(n)=k}} 1 \nonumber
\end{eqnarray}
in the form of a series similar to that of (2.5). \\ \\
If $k$ is an arbitrary but fixed positive integer, then it can be shown by induction on $k$ starting with (2.5), that there exist polynomials $Q_{j,k}(X)$, of degree $\leq k-1$, with $Q_{1,k}(X)$ of degree $k-1$ and leading coefficient 1, such that 
\begin{multline}
    \pi_k(x) = \frac{xQ_{1,k}(\log\log x)}{\log x}+\frac{xQ_{2,k(\log\log x)}}{\log^2x}+\cdots  +\frac{xQ_{\nu,k}(\log\log x)}{\log^{\nu}x}+O_{\nu}\left(\frac{x(\log\log x)^k}{\log^{\nu+1}x}\right).
\end{multline}
\textbf{Remark:} For $k$, an arbitrary but fixed positive integer, Landau obtained by induction on $k$, the asymptotic estimate
\begin{eqnarray}
    \pi_k(x) \sim \frac{x(\log\log x)^{k-1}}{\log x}
\end{eqnarray}
as a consequence of the Prime Number Theorem. If $k$ varies with $x$, then it becomes cumbersome to keep track of the error terms in the inductive argument, yet Sathe (see \cite{SatheI}, \cite{SatheII}) skillfully managed to obtain asymptotic uniform estimates for $\pi_k(x)$ for $k \leq B\log\log x$, where $B$ is an arbitrary constant. Selberg \cite{Sel54} subsequently showed how to get Sathe's results painlessly by an analytic method. Delange later improved upon Selberg's method and the Selberg-Delange analytic method \cite{Tbook} yields a series representation similar to (2.6) for $k \leq B\log\log x$. Our emphasis here is that, if $k$ is arbitrary but fixed, then (2.6) can be derived elementarily by induction on $k$ using the strong form of the Prime Number Theorem as given in (2.5). \\ \\
In the series representation for $M_{\omega^k}(x)$ and $m_{\omega^k}(x)$ given below, terms involving $\log\log x$ appear only for $k \geq 2$. Thus, the analysis for case $k=2$, although bearing similarity to the case $k=1$ in Alladi-Johnson \cite{AJ24}, is different in a crucial way. The analysis of all cases $k \geq 2$ are all quite similar, but the details get more complicated as $k$ increases. For these reasons, we shall present here all the details for the case $k=2$, and state the corresponding results for $k\geq 3$ sketching only the main ideas in the proof. \\ \\
We begin by noting that if $\omega(n)=r \geq 1$, then
\begin{eqnarray}
    \sum_{d|n} \mu(d)z^{\omega(d)} = (1-z)^r.
\end{eqnarray}
So
\begin{eqnarray}
    \sum_{d|n} \mu(d)\omega(d)z^{\omega(d)} = z\frac{d}{dz}\left(\sum_{d|n}\mu(d)z^{\omega(d)}\right) = -zr(1-z)^{r-1}.
\end{eqnarray}
From (2.9), we see that 
\begin{eqnarray}
    \sum_{d|n} \mu(d)\omega^2(d)z^{\omega(d)-1} = \frac{d}{dz}\{-zr(1-z)^{r-1}\} = -r(1-z)^{r-1}+zr(r-1)(1-z)^{r-2}.
\end{eqnarray}
Thus,
\begin{eqnarray}
    \sum_{d|n}\mu(d)\omega^2(d) = \left(\sum_{d|n}\mu(d)\omega^2(d)z^{\omega(d)-1}\right)\Bigg|_{z=1} = -\chi_1(n)+2\chi_2(n),
\end{eqnarray}
where $\chi_1$ and $\chi_2$ are the characteristic functions of integers for which $\omega(n)=1$ and $\omega(n)=2$ respectively. Therefore, by M\"obius inversion, we have
\begin{eqnarray}
    \mu(n)\omega(n)^2 = \sum_{d|n}\left(-\chi_1(d)+2\chi_2(d)\right)\mu\left(\frac{n}{d}\right).
\end{eqnarray}
In what follows, (2.12) will be used to estimate $M_{\omega^2}(x)$ asymptotically using the hyperbola method. In order to do this, we need some results on the M\"obius function. \\ \\
The Dirichlet series associated with the M\"obius function is
\begin{eqnarray}
    \sum_{n=1}^{\infty}\frac{\mu(n)}{n}=\frac{1}{\zeta(s)}, \quad \text{for $\sigma>1$},
\end{eqnarray}
where $\zeta(s)$ is the Riemann zeta function. By using bounds for $\frac{1}{\zeta(s)}$, the standard analytic approach to obtain the strong form of the Prime Number Theorem can be employed to derive the following bound for $M(x)$ - see for instance, Tenenbaum (\cite{Tbook}, p. 217):
\begin{eqnarray}
    M(x) = \sum_{n\leq x} \mu(n)\ll xe^{-c_1\sqrt{\log x}}.
\end{eqnarray}
The same method also yields
\begin{eqnarray}
    m(x) = \sum_{n\leq x}\frac{\mu(n)}{n} \ll e^{-c_1\sqrt{\log x}},
\end{eqnarray}
and these are quantitative forms of (1.1) and (1.2). In (2.14) and (2.15), we have used the same constant $c_1$, because if we had two different positive constants, we could choose the minimum of these two as $c_1$ for (2.14) and (2.15), we have for each positive integer j,
\begin{eqnarray}
    \sum_{n=1}^{\infty} \frac{\mu(n)\log^jn}{n} = \lambda_j
\end{eqnarray}
is convergent, and its quantitative form
\begin{eqnarray}
    \sum_{n\leq x} \frac{\mu(n)\log^jn}{n} = \lambda_j+O(e^{-c_1\sqrt{\log x}}).
\end{eqnarray}
The constant $c_1=c_1(j)$ in the exponential of (2.17) depends on $j$, but in what follows, we will use (2.17) only for $j \leq \nu$, where $\nu$ will be an arbitrary but fixed positive integer. So, we will use
\begin{eqnarray}
    c_1:=\min\{c_1(0),c_1(1),\cdots  ,c_1(\nu)\}.
\end{eqnarray}
Note that (2.17) can be deduced from the method that yields (2.14) or by partial summation using (2.14). It is known that
\begin{eqnarray}
    \sum_{n=1}^{\infty} \frac{\mu(n)\log n}{n} = \lambda_1 = -1.
\end{eqnarray}
This is the only value besides $\lambda_0=0$ that we need. The actual values of $\lambda_j$ for $j \geq 2$, which can be written in terms of the successive derivative values of $\frac{1}{\zeta(s)}$ at $s=1$, will not concern us. We are now ready to prove the following theorem: \\ \\
\textbf{Theorem 2.1:} \textit{Let $\nu$ be an arbitrary but fixed integer. Then there exist linear polynomials $Q_{j,2}^*(X)$ such that}
\begin{eqnarray}
    M_{\omega^2}(X) = -\frac{x(2\log\log x+c_2)}{\log^2x}+\sum_{j=3}^{\nu} \frac{xQ_{j,2}^*(\log\log x)}{\log^jx}+O_{\nu}\left(\frac{x\log\log^{2\nu+4}x}{\log^{\nu+1}x}\right). \nonumber
\end{eqnarray}
\textbf{Proof:}  We shall adapt the method in Alladi-Johnson \cite{AJ24} to prove \textit{Theorem 2.1}. Begin by using (2.12) to write
\begin{eqnarray}
    M_{\omega^2}(x) = \sum_{n\leq x}\mu(n)\omega(n)^2 = \sum_{n\leq x}\sum_{d|n}\mu\left(\frac{n}{d}\right) (2\chi_2(d)-\chi_1(d)).
\end{eqnarray}
Using the hyperbola method, we break up the double sum in (2.20) into
\begin{align}
     M_{\omega^2}(x)  &= \sum_{m\leq T}\mu(m)\sum_{h \leq \frac{x}{m}} (2\chi_2(h)-\chi_1(h)) + \sum_{h\leq \frac{x}{T}}(2\chi_2(h)-\chi_1(h))\sum_{T<m\leq \frac{x}{h}}\mu(m) \nonumber \\
     &= \Sigma_1+\Sigma_2,
\end{align}
where $T$ will be chosen optimally below so that
\begin{eqnarray}
    \log T = o(\log x).
\end{eqnarray}
We can easily bound $\Sigma_2$ using (2.14) and the monotone increasing property of 
\begin{eqnarray}
    R_1(x) = xe^{-c_1\sqrt{\log x}},\quad\text{for $x\geq x_0$}. \nonumber
\end{eqnarray}
That is, (2.14) yields
\begin{eqnarray}
    \Sigma_2 \ll \sum_{h \leq \frac{x}{T}}\frac{x(2\chi_2(h)-\chi_1(h))}{he^{c_1\sqrt{\log \left(\frac{x}{h}\right)}}} \ll \frac{x\log\log x}{e^{c_1\sqrt{\log T}}}.
\end{eqnarray}
The estimation of $\Sigma_1$ is more involved. We note that (2.5) and (2.6) yield
\begin{align}
    \sum_{h\leq x} \{2\chi_2(h)-\chi_1(h)\} &=2\pi_2(x)-\pi(x)+O(\sqrt{x}) \nonumber \\ 
    &= \frac{x(2\log\log x+c_2)}{\log x}+\frac{xF_{2,2}(\log\log x)}{\log^{2}x}+\cdots  +\frac{xF_{\nu,2}(\log\log x)}{\log^{\nu}x}+O_{\nu}\left(\frac{x\log\log x}{\log^{\nu+1}x}\right),
\end{align}
where, for $2 \leq j \leq \nu$,
\begin{eqnarray}
    F_{j,2}(\log\log x) = 2Q_{j,2}(\log\log x)-(j-1)!. \nonumber
\end{eqnarray}
So, from (2.21) and (2.22) we get 
\begin{multline}
    \Sigma_1 = \sum_{m \leq T}\mu(m)\left(\frac{x\left(2\log\log \left(\frac{x}{m}\right)+c_2\right)}{m\log\left(\frac{x}{m}\right)}\right) +\sum_{m\leq T}\left(\sum_{j=2}^{\nu}\frac{\mu(m)xF_{j,2}\left(\log\log \left(\frac{x}{m}\right)\right)}{m\log^j\left(\frac{x}{m}\right)}\right)+O_{\nu}\left(\frac{x\log\log x\log T}{\log^{\nu +1}\left(\frac{x}{T}\right)}\right).
\end{multline}
Denote by
\begin{eqnarray}
    \Sigma_{1,j} = x\sum_{m\leq T}\frac{\mu(m)F_{j,2}\left(\log\log \left(\frac{x}{m}\right)\right)}{m\log^j\left(\frac{x}{m}\right)}, \quad \text{for $j=1,2,\cdots  ,\nu$}, 
\end{eqnarray}
where $F_{1,2}(\log\log x)=2\log\log x+c_2$. \\ \\
We first estimate $\Sigma_{1,1}$. To this end, we note that since $m \leq T$ and (2.22) holds, we have 
\begin{align}
    \log\log \left(\frac{x}{m}\right) &= \log(\log x-\log m) \nonumber \\
    &= \log\log x+\log\left(1-\frac{\log m}{\log x}\right) \nonumber \\
    &=\log\log x - \sum_{j=1}^{\nu}\frac{\log^j m}{j\log^j x}+O\left(\frac{\log^{\nu+1}m}{\log^{\nu+1}x}\right),
\end{align}
and
\begin{align}
    \frac{1}{\log\left(\frac{x}{m}\right)} &= \frac{1}{\log x\left(1-\frac{\log m}{\log x}\right)} \nonumber \\
    &= \sum_{j=0}^{\nu} \frac{\log^j m}{\log^{j+1}x}+O\left(\frac{\log^{\nu+1}m}{\log^{\nu+2}x}\right).
\end{align}
So, from (2.27) and (2.28), we get
\begin{eqnarray}
    \frac{\log\log \left(\frac{x}{m}\right)}{\log \left(\frac{x}{m}\right)} = \frac{\log\log x}{\log x} + \sum_{j=1}^{\nu}\frac{(\log\log x-H_j)\log^j m}{\log^{j+1}x} +O_{\nu} \left(\frac{\log\log x\log^{\nu+1}}{\log^{\nu+2}x}\right),
\end{eqnarray}
where the $H_j$ are the Harmonic numbers given by
\begin{eqnarray}
    H_j=\sum_{i=1}^j \frac{1}{i}. \nonumber
\end{eqnarray}
Thus, from (2.26) and (2.29), we see that
\begin{align}
    \Sigma_{1,1} &= \frac{x(2\log\log x+c_2)}{\log x}\sum_{m \leq T}\frac{\mu(m)}{m} +\frac{x(2\log\log x+c_2-H_1)}{\log^2x}\sum_{m\leq T}\frac{\mu(m)\log m}{m}\nonumber \\
    &\hspace{1.3cm}+\sum_{j=2}^{\nu}\frac{x(\log\log x+c_2-H_j)}{\log^{j+1}x}\sum_{m \leq T}\frac{\mu(m)\log^jm}{m} + O_{\nu} \left(\frac{x\log\log x\log^{\nu+2}T}{\log^{\nu+2}x}\right) \nonumber \\
    &= -\frac{x(2\log\log x+c_2-H_1)}{\log^2x}+\sum_{j=2}^{\nu}\frac{x(\log\log x+c_2-H_j)}{\log^{j+1}x} \nonumber \\
    &\hspace{1.3cm}+O_{\nu}\left(\frac{x\log\log x\log^{\nu+2}T}{\log^{\nu+2}x}\right)+O_{\nu}\left(\frac{x\log\log x}{e^{c_1\sqrt{\log x}}}\right) ,
\end{align}
in view of (2.15), (2.17) and (2.19). \\ \\
Note that the leading term in $\Sigma_{1,1}$ is $\frac{x(\log\log x+c_2)}{\log^2x}$. We can treat the sums $\Sigma_{1,j}$ for $j \geq 2$ in a similar fashion and obtain a series in decreasing powers of $\log x$ with leading term
\begin{eqnarray}
    \frac{F_{j,2}(\log\log x)}{\log^{j+1}x} \nonumber
\end{eqnarray}
and error terms as in (2.30), because $\lambda_0=0$, where $F_{j,2}(X)$ is a linear function of $X$. Finally, if we sum all such expressions for $\Sigma_{1,j}$, we will get
\begin{eqnarray}
    \Sigma_1 = -\frac{x(2\log\log x+c_2)}{\log^2x} + \sum_{j=3}^{\nu}\frac{xQ^*_{j,2}(\log\log x)}{\log^jx}+O_{\nu}\left(\frac{x\log\log^{2\nu+4}x}{\log^{\nu+1}x}\right)+O_{\nu}\left(\frac{x\log\log x}{e^{c_1\sqrt{\log T}}}\right).
\end{eqnarray}
At this stage, we choose 
\begin{eqnarray}
    T=e^{U^2(\log\log x)^2}, \quad \text{with $c_1U=\nu+1$}.
\end{eqnarray}
With this choice of $T$, \textit{Theorem 2.1} follows from (2.21), (2.23) and (2.31). \qed \\ \\
Next, for the purpose of proving \textit{Theorem 2.2}, and also for subsequent use in this paper, and for use in future related work, we will use \textit{Theorem A} below due to Alladi-Johnson \cite{AJ24}, which is a variant of Axer's theorem (see Hardy \cite{HarDiv}; p. 378) for function which need not be bounded: \\ \\
\textbf{Theorem A:} \textit{Let $\{a_n\}_{n=1}^{\infty}$ be a sequence of reals such that}
\begin{eqnarray}
    A(x):=\sum_{n \leq x}a_n \ll x\eta(x),
\end{eqnarray}
\textit{where}
\begin{equation}
   \tag{2.34a} \eta(x) \to 0 \quad \text{as} \quad x\to \infty,
\end{equation}
\textit{and}
\begin{equation}
    x\eta(x) \quad \text{is an increasing function $\to \infty$} \quad \text{as} \quad x\to \infty. \tag{2.34b}
\end{equation}
\textit{Suppose also that}
\begin{equation}
    \sum_{n \leq x}|a_n| \ll x\beta(x), \tag{2.35}
\end{equation}
\textit{where}
\begin{equation}
    \beta(x) \quad \text{is an increasing function,} \quad \text{and} \quad \beta(x)\eta(x) \to 0 \quad \text{as}\quad x\to \infty. \tag{2.36}
\end{equation}
\textit{Then}
\begin{equation}
    \sum_{n\leq x}a_n\left\{\frac{x}{n}\right\} \ll x\sqrt{\eta(x)\beta(x)} = o(x), \tag{2.37}
\end{equation}
\textit{where $\{x\}$ is the fractional part of $x$}. \\ \\
\textbf{Remark:} For a proof, see \cite{AJ24}. Axer \cite{axer} stated his theorem for bounded functions (sequences) $a_n$, but the method can be applied with $|a_n|$ having at most a slowly growing unbounded average as given by (2.35) and (2.36). In the version of Axer's theorem in Hardy (\cite{HarDiv}, p. 378), the function $\{x\}$ is replaced by a more general function $\chi(x)$ of bounded variation on finite intervals, but $|a_n|$ is assumed to have a bounded average. \textit{Theorem A} can be generalized by replacing $\{x\}$ with such a function $\chi(x)$, but the version of \textit{Theorem A} given above suffices for our purpose here and in our subsequent work. \\ \\
\textbf{Theorem 2.2:} \textit{With $\{w\}$ denoting the fractional part of $w$, we have}
\begin{eqnarray}
    \sum_{n \leq x}\mu(n)\omega(n)^2 \left\{\frac{x}{n}\right\} \ll \frac{x(\log\log x)^{\frac{3}{2}}}{\log x}. \nonumber
\end{eqnarray}
\textbf{Proof:} To prove \textit{Theorem 2.2}, choose $a_n=\mu(n)\omega(n)^2$ in \textit{Theorem A}. By \textit{Theorem 2.1}, we can take
\begin{eqnarray}
    \eta(x) = \frac{\log\log x}{\log^2x} \nonumber
\end{eqnarray}
in (2.33). Thus, (2.34a) and (2.34b) are satisfied. In view of the well known estimate
\begin{eqnarray}
    \sum_{1\leq n \leq x}\omega(n)^2 \ll x(\log\log x)^2, \nonumber
\end{eqnarray}
we can take
\begin{eqnarray}
    \beta(x) = (\log\log x)^2 \nonumber
\end{eqnarray}
in (2.35), since $|a_n|\leq \omega(n)^2$. Thus, (2.36) is satisfied. \textit{Theorem 2} then follows from (2.37) of \textit{Theorem A}. \qed \\ \\
Since, \textit{Theorem 2.2} deals with the fractional part function as the weight, we establish next the corresponding result with the weight as the integral part function: \\ \\
\textbf{Theorem 2.3:} \textit{Let $[w]$ denote the integral part of $w$. Then for each positive integer $\nu$, we have}
\begin{align}
    \sum_{n\leq x}\mu(n)\omega(n)^2\left[\frac{x}{n}\right] = \frac{x(Q_{1,2}(\log\log x)-1)}{\log x} &+\frac{x(Q_{2,2}(\log\log x)-1)}{\log^2x}+\frac{xQ_{2,3}(\log\log x)-2!}{\log^3x}+\cdots  \nonumber \\
    &+\frac{x(Q_{\nu,2}(\log\log x)-(\nu-1)!)}{\log^{\nu}x}+O\left(\nu!\cdot\frac{x(\log\log x)^{\nu+1}}{\log^{\nu+1}x}\right). \nonumber
\end{align}
\textbf{Proof:} Note that (2.11) yields
\begin{eqnarray}
    \sum_{n\leq x}\mu(n)\omega(n)^2\left[\frac{x}{n}\right] = \sum_{n\leq x}\sum_{d|n} \mu(d)\omega^2(d) = \sum_{n\leq x}(2\chi_2(n)-\chi_1(n)), \nonumber
\end{eqnarray}
and from this, \textit{Theorem 2.3} follows using (2.5) and (2.6). We note here that the $Q_{j,2}(\log\log x)$ are linear polynomials in $\log\log x$. \qed \\ \\
We now prove the main result of this section which we state in two parts. The reason for this split will be clear at the end of the proof. \\ \\
\textbf{Theorem 2.4:}
\begin{itemize}
    \item[(i)] \textit{We have}
    \begin{eqnarray}
        m_{\omega^2}(x):=\sum_{n\leq x}\frac{\mu(n)\omega(n)^2}{n} = O\left(\frac{(\log\log x)^{\frac{3}{2}}}{\log x}\right) .\nonumber
    \end{eqnarray}
    \textit{Consequently}
    \begin{eqnarray}
        \sum_{n=1}^{\infty} \frac{\mu(n)\omega(n)^2}{n} = \sum_{n=2}^{\infty} \frac{\mu(n)\omega(n)^2}{n} = 0 .\nonumber
    \end{eqnarray}
    \item[(ii)] \textit{More precisely, there exist linear polynomials $F_{j,2}^*(X)$, such that, for each positive integer $\nu$, we have}
    \begin{eqnarray}
        m_{\omega^2}(x) = \frac{2\log\log x+c_2+2}{\log x}+\frac{F_{2,2}^*(\log\log x)}{\log^2x}+\cdots  +\frac{F^*_{\nu,2}(\log\log x)}{\log^{\nu}x} + O_{\nu}\left(\frac{\log\log^{2\nu+6}x}{\log^{\nu+1}x}\right) .\nonumber
    \end{eqnarray}
\end{itemize}
\textbf{Proof of (i):} From \textit{Theorem 2.2} and \textit{2.3}, we get
\begin{align}
    \sum_{n \leq x}\mu(n)\omega(n)^2\frac{x}{n} &= \sum_{n\leq x}\mu(n)\omega(n)^2\left[\frac{x}{n}\right] + \sum_{n\leq x}\mu(n)\omega(n)^2\left\{\frac{x}{n}\right\} \nonumber \\
    &\ll \frac{x(\log\log x)^{\frac{3}{2}}}{\log x}+\frac{x\log\log x}{\log x} \ll \frac{x(\log\log x)^{\frac{3}{2}}}{\log x} .\tag{2.38}
\end{align}
By canceling $x$ on both extremes of (2.38), we get the bound for $m_{\omega^2}(x)$ in \textit{Theorem 2.4 (i)}. By letting $x \to \infty$ in this bound, we get the second assertion of \textit{Theorem 2.4 (i)}. \qed \\ \\
\textbf{Proof of (ii):} We start with the representation 
\begin{equation}
    m_{\omega^2}(x) =\sum_{n\leq x}\frac{\mu(n)\omega(n)^2}{n} = \int_1^x \frac{dM_{\omega^2}(t)}{t}. \tag{2.39}
\end{equation}
Note that $M_{\omega^2}(t)=m_{\omega^2}(t)=0$ for $t<2$. Integration-by-parts of the Stieltjes integral in (2.39) gives
\begin{align}
    m_{\omega^2}(t) &= \frac{M_{\omega^2}(t)}{t}\Big|_{1}^x + \int_1^x \frac{M_{\omega^2}(t)}{t^2}dt \nonumber \\
    &= \frac{M_{\omega^2}(x)}{x}+ \int_1^x \frac{M_{\omega^2}(t)}{t^2}dt. \tag{2.40}
\end{align}
We know from \textit{Theorem 2.4 (i)} that $m_{\omega^2}(\infty)=0$, and from \textit{Theorem 2.1} that $M_{\omega^2}(x)=o(x)$. So, by letting $x\to \infty$ in (2.40), we deduce that
\begin{equation}
    \int_1^{\infty} \frac{M_{\omega^2}(t)}{t^2}dt = 0. \tag{2.41}
\end{equation}
In view of (2.41), we may rewrite (2.40) as
\begin{equation}
    m_{\omega^2}(x) = \frac{M_{\omega^2}(x)}{x}-\int_x^{\infty} \frac{M_{\omega^2}(t)}{t^2}dt. \tag{2.42}
\end{equation}
From \textit{Theorem 2.1}, we see that 
\begin{equation}
    \frac{M_{\omega^2}(x)}{x} = -\frac{2\log\log x+c_2}{\log^2x} + \sum_{j=3}^{\nu}\frac{Q_{j,2}^*(\log\log x)}{\log^jx}+O_{\nu}\left(\frac{\log\log^{2\nu+4}x}{\log^{\nu+1}x}\right). \tag{2.43}
\end{equation}
As for the integral on the right in (2.42), if we substitute the expression in \textit{Theorem 2.1} for $M_{\omega^2}(t)$ in the integral with $\nu$ replaced $\nu+1$, we will get
\begin{equation}
    -\int_x^{\infty} \frac{M_{\omega^2}(t)}{t^2}dt = \int_{x}^{\infty}\frac{(2\log\log t+c_2)}{t\log^2t}dt- \sum_{j=3}^{\nu+1}\int_x^{\infty}\frac{Q_{j,2}^*(\log\log t)}{t\log^jt}dt+O_{\nu}\left(\int_x^{\infty}\frac{\log\log^{2\nu+4}t}{t\log^{\nu+1}t}dt\right). \tag{2.44}
\end{equation}
In (2.44), the first integral on the right is what gives the leading term in the series expansion for $m_{\omega^2}(x)$, and it is
\begin{align}
    \int_{x}^{\infty} \frac{2(\log\log t+c_2)}{t\log^2t}dt &= \int_{x}^{\infty} \frac{2\log\log t}{t\log^2t}dt+2c_2\int_{x}^{\infty} \frac{1}{t\log^2t}dt \nonumber \\
    &= \int_{x}^{\infty} \frac{2\log\log t}{t\log^2t}dt+ \frac{c_3}{\log x}, \tag{2.45}
\end{align}
where $c_3=2c_2$, another constant. For the integral on the right of (2.45), we use the substitution $v=\log\log t$ to get
\begin{equation}
    \int_{x}^{\infty} \frac{2\log\log t}{t\log^2t}dt = 2\int_{\log\log x}^{\infty} v e^{-v}dv = 2(-v e^{-v}-e^{-v})\big|_{\log\log x}^{\infty} = \frac{2\log\log x+2}{\log x}. \tag{2.46}
\end{equation}
So, from (2.45) and (2.46), we get
\begin{equation}
    \int_x^{\infty} \frac{(2\log\log t+c_2)dt}{t\log^2t} = \frac{2\log\log x+c_3+2}{\log x} .\tag{2.47}
\end{equation}
The other integrals can be worked out in an analogous fashion. If we denote the generic linear polynomial $Q_{j,2}^*(\log\log x)$ by $a_j\log\log x+b_j$, then we get
\begin{align}
    \int_{x}^{\infty}\frac{Q_{j,2}^*(\log\log t)}{\log^jt}dt &= \int_x^{\infty}\frac{a_j\log\log t}{t\log^jt}dt+b_j\int_x^{\infty}\frac{dt}{t\log^jt} \nonumber \\
    &= a_j\int_{x}^{\infty}ve^{-(j-1)v}dv+\frac{b_j}{(j-1)\log^jx} \nonumber \\
    &= \frac{a_j\log\log x}{(j-1)\log^{j-1}x}+\frac{a_j}{(j-1)^2\log^{j-1}x}+\frac{b_j}{(j-1)\log^{j-1}x}, \tag{2.48}
\end{align}
using the substitution $\log\log t = v$ in the integral involving $a_j$, and then integrating-by-parts. The terms in (2.48) have to be summed from $j=3$ to $j=\nu+1$, and then combined with the terms in (2.44). Finally,
\begin{equation}
    O_{\nu}\left(\int_x^{\infty}\frac{\log\log^{2\nu+4}t}{t\log^{\nu+1}t}dt\right) \ll \frac{\log\log^{2\nu+6}x}{\log^{\nu+1}x}. \tag{2.49}
\end{equation}
Thus, \textit{Theorem 2.4 (ii)} follows from the equations (2.44) to (2.49). \qed 
\subsection{Results involving $\mu(n)\omega(n)^k$ for $k \geq 3$}
We now state the extensions of \textit{Theorems 2.1, 2.2, 2.3} and \textit{2.4,} for sums involving $\mu(n)\omega(n)^k$ when $k \geq 3$. We title these results below as \textit{Theorems 2.1-k, 2.2-k, 2.3-k} and \textit{2.4-k}, respectively. We provide the main ideas needed for establishing these four theorems for $k\geq 3$, but do not go through the details of the proof because of the similarity with the case $k=2$ discussed above. \\ \\
In order to prove \textit{Theorem 2.1}, we started with (2.8) and considered $z\frac{d}{dz}$ of the expressions in (2.8) to get (2.9). We then considered $\frac{d}{dz}$ of the expressions in (2.9) to get (2.10) and evaluated this at $z=1$ to get identity (2.11). It is with that we started the proof of \textit{Theorem 2.1} (see (2.20)). In order to estimate 
\begin{eqnarray}
    M_{\omega^k}(x):=\sum_{n\leq x} \mu(n)\omega(n)^k, \nonumber
\end{eqnarray}
for $k\geq 3$, we need an expression analogous to (2.11) for the sum
\begin{eqnarray}
    \sum_{d|n} \mu(d)\omega^k(d). \nonumber
\end{eqnarray}
To get this, we need to multiply the expressions in (2.10) by $z$ and then consider $\frac{d}{dz}$ of that expression, and iterate this process of multiplying by $z$ and considering $\frac{d}{dz}$ until we get
\begin{equation}
    \sum_{d|n} \mu(d)\omega^k(d)z^{\omega(d)-1} = f_k(z), \tag{2.50}
\end{equation}
where $f_k(z)$ is the analogue of the expression on the right hand side of (2.10) for $k \geq 3$. Then, by evaluating $f_k(1)$, we end up with the identity
\begin{equation}
    \sum_{d|n} \mu(d)\omega^k(d) = \sum_{j=1}^{k} \delta_{j,k}\chi_j(n), \tag{2.51}
\end{equation}
where $\delta_{j,k}$ are integer constants and $\chi_j$ is the characteristic function of integers $m$ for which $\omega(m)=j$. A few calculations give us the following formula for the constants $\delta_{j,k}$: for  a fixed positive integer $k$ and a positive integer $j$ satisfying $1 \leq j \leq k$ 
\begin{eqnarray}
    \delta_{j,k}:=\begin{cases}
        -1,\quad \text{for $j=1$} \\
        (-1)^kk!, \quad \text{for $j=k$} \\
        j(\delta_{j,k-1}-\delta_{j-1,k-1}), \quad \text{for $1<j<k$.} 
    \end{cases} \nonumber
\end{eqnarray}
From (2.51) and M\"obius inversion, we get
\begin{equation}
    M_{\omega^k}(x) = \sum_{n\leq x}\mu(n)\omega(n)^k = \sum_{n\leq x}\sum_{d|n}\left(\sum_{j=1}^k\delta_{j,k}\chi_j(d)\mu\left(\frac{n}{d}\right)\right). \tag{2.52}
\end{equation}
The idea is to estimate $M_{\omega^k}(x)$ by the hyperbola method, appealing to the bound for $M(x)$ given in (2.14) and using the fact that
\begin{eqnarray}
    \sum_{n\leq x} \chi_j(n) = \pi_j(x)+O(\sqrt{x}), \nonumber
\end{eqnarray}
together with the evaluation of $\pi_j(x)$ given by (2.6). Here, instead of (2.27), we need to use
\begin{align}
    \left\{\log\log \left(\frac{x}{m}\right)\right\}^t =\left\{\log\log x - \sum_{j=1}^{\nu}\frac{\log^j m}{j\log^j x}+O\left(\frac{\log^{\nu+1}m}{\log^{\nu+1}x}\right)\right\}^t,\quad \text{for $t=1,2,3,\cdots  ,k-1$}, \nonumber
\end{align}
and follow the argument in the proof of \textit{Theorem 2.1}.
The leading term, when expanding the right hand side of the above expression, will be $(\log\log x)^t$.
The splitting of the sum in (2.52) into $\Sigma_{1,k}$ and $\Sigma_{2,k}$ is identical to the split into $\Sigma_1$ and $\Sigma_2$ in (2.21), and the parameter $T$ will again be chosen to satisfy (2.22). Indeed, the choice of $T$ as in (2.32) works for all $k \geq 3$, and the result that we get is the following: \\ \\
\textbf{Theorem 2.1-k:} \textit{Let $k \geq3$, and $\nu$ be an arbitrary fixed integer. Then there exist polynomials $Q_{j,k}^*(X)$ of degree $\leq  k-1$ in $X$ with $Q_{2,k}^*(X)$ of degree $k-1$ such that}
\begin{eqnarray}
    M_{\omega^k}(x) = \frac{xQ_{2,k}^*(\log\log x)}{\log^2x}+\sum_{j=3}^{\nu}\frac{xQ_{j,k}^*(\log\log x)}{\log^jx}+O_{\nu}\left(\frac{x(\log\log x)^{k\nu+2k}}{\log^{\nu+1}x}\right). \nonumber
\end{eqnarray}
From \textit{Theorem 2.1-k} and \textit{Theorem A}, we get: \\ \\
\textbf{Theorem 2.2-k:} \textit{With $\{w\}$, we have for $k\geq 3$}
\begin{eqnarray}
    \sum_{n\leq x}\mu(n)\omega(n)^k\left\{\frac{x}{n}\right\} \ll_k \frac{x(\log\log x)^{\frac{2k-1}{2}}}{\log x}. \nonumber
\end{eqnarray}
\textbf{Proof:} To prove \textit{Theorem 2.2-k}, choose $a_n = \mu(n)\omega(n)^k$ in \textit{Theorem A}. By \textit{Theorem 2.1-k}, we know that in this case, we can take
\begin{eqnarray}
    \eta(x) = \frac{(\log\log x)^{k-1}}{\log^2 x} \nonumber
\end{eqnarray}
in (2.33). Thus, (2.34a) and (2.34b) are satisfied. In view of the well known estimate 
\begin{eqnarray}
    \sum_{1\leq n\leq x} \omega(n)^k \ll x(\log\log x)^k, \nonumber
\end{eqnarray}
we can take 
\begin{eqnarray}
    \beta(x) = (\log\log x)^k \nonumber
\end{eqnarray}
in (2.35), since $|a_n| \leq \omega(n)$. Thus, (2.36) is satisfied. \textit{Theorem 2.2-k} then follows from (2.37) of \textit{Theorem A}. \qed \\ \\
Next, we state the corresponding result with the weight as the integral part function. \\ \\
\textbf{Theorem 2.3-k:} \textit{Let $[w]$ denote the integral part of $w$, and $k \geq 3$. Then we have}
\begin{eqnarray}
    \sum_{n \leq x}\mu(n)\omega(n)^k \left[\frac{x}{n}\right] \sim \frac{\delta_{k,k}(\log\log x)^{k-1}}{\log x}.\nonumber
\end{eqnarray}
\textbf{Proof:} By (2.51), we have
\begin{eqnarray}
    \sum_{n\leq x}\mu(n)\omega(n)^k\left[\frac{x}{n}\right] = \sum_{n\leq x}\sum_{d|n}\mu(d)\omega^k(d) = \sum_{n\leq x}\sum_{j=1}^{k}\delta_{j,k}\chi_j(n)\sim \frac{\delta_{k,k}(\log\log x)^{k-1}}{\log x} ,\nonumber
\end{eqnarray}
in view of (2.6). \qed\\ \\
\textbf{Remark:} Like \textit{Theorem 2.3}, there is a series expansion for the left hand side of \textit{Theorem 2.3-k}, but this involves writing the expansions for each of
\begin{eqnarray}
    \sum_{n\leq x}\chi_j(n), \quad \text{for $j=1,2,\cdots  ,k$}, \nonumber
\end{eqnarray}
and adding them. Such a refinement is not needed for our purpose. \\ \\
We now state the main result (\textit{Theorem 2.4-k}) of this section for the general case $k \geq 3$. We state the theorem in two parts just as we did for \textit{Theorem 2.4}. Since \textit{Theorem 2.4-k} can be proved by a suitable extension of the method used to prove \textit{Theorem 2.4}, we only provide the main ideas needed to prove \textit{Theorem 2.4-k} and skip the details.\\ \\
\textbf{Theorem 2.4-k:} \textit{Let $k \geq 3$ and $\nu$ an arbitrary but fixed positive integer.}
\begin{itemize}
    \item[(i)] \textit{We have}
    \begin{eqnarray}
        m_{\omega^k}(x):= \sum_{n\leq x}\frac{\mu(n)\omega(n)^k}{n} = O\left(\frac{(\log\log x)^{\frac{2k-1}{2}}}{\log x}\right). \nonumber
    \end{eqnarray}
    \textit{Consequently}
    \begin{eqnarray}
        \sum_{n=1}^{\infty}\frac{\mu(n)\omega(n)^k}{n} = \sum_{n=2}^{\infty}\frac{\mu(n)\omega(n)^k}{n} = 0. \nonumber
    \end{eqnarray}
    \item[(ii)] \textit{More precisely, there exist polynomials $F_{j,k}^*(X)$ of degree $\leq k-1$ in $X$, with $F^*_{1,k}$ being of degree $k-1$, such that for each positive integer $\nu$, we have}
    \begin{eqnarray}
        m_{\omega^k}(x) = \frac{F_{1,k}^*(\log\log x)}{\log x}+\frac{F_{2,k}^*(\log\log x)}{\log^2 x}+\cdots  +\frac{F_{\nu,k}^*(\log\log x)}{\log^{\nu} x}+O_{\nu}\left(\frac{\log\log^{k\nu+3k}x}{\log^{\nu+1}x}\right). \nonumber
    \end{eqnarray}
\end{itemize}
\textbf{Proof of Theorem 2.4-k (i):} From \textit{Theorem 2.2-k} and \textit{3-k}, we get 
\begin{align}
    \sum_{n\leq x}\mu(n)\omega(n)^k\frac{x}{n} &=\sum_{n\leq x } \mu(n)\omega(n)^k\left[\frac{x}{n}\right] +\sum_{n\leq x}\mu(n)\omega(n)^k\left\{\frac{x}{n}\right\} \nonumber \\
    &\ll \frac{x(\log\log x)^{k-1}}{\log x}+\frac{x(\log\log x)^{\frac{2k-1}{2}}}{\log x} \nonumber \\
    &\ll \frac{x(\log\log x)^{\frac{2k-1}{2}}}{\log x}. \tag{2.53}
\end{align}
By canceling $x$ on both extremes of (2.53), we get the bound for $m_{\omega^k}(x)$ in \textit{Theorem 2.4-k (i)}. By letting $x\to \infty$ in this bound, we get the second assertion in \textit{Theorem 2.4-k (i)}. \qed \\ \\
\textbf{Proof of Theorem 2.4-k (ii):} We start with the representation
\begin{equation}
    m_{\omega^k}(x) = \sum_{n\leq x} \frac{\mu(n)\omega(n)^k}{n} = \int_1^{\delta}\frac{dM_{\omega^k}(t)}{t} .\tag{2.54}
\end{equation}
Note that $M_{\omega^k}(t)=m_{\omega^k}(t)=0$ for $t<2$. integration-by-parts pf the Stieltjes integral in (2.54) gives 
\begin{equation}
    m_{\omega^k}(x) = \frac{M_{\omega^k}(t)}{t}\Big|_1^x+\int_1^x\frac{M_{\omega^k}(t)}{t^2}dt = \frac{M_{\omega^k}(x)}{x}+\int_1^x\frac{M_{\omega^k}(t)}{t^2}dt. \tag{2.55}
\end{equation}
We know from \textit{Theorem 2.4-k (i)} that $m_{\omega^k}(\infty)=0$, and from \textit{Theorem 2.1-k} that $M_{\omega^k}(x)=o(x)$. So, by letting $x \to \infty$ in (2.55), we deduce that
\begin{equation}
    \int_1^{\infty} \frac{M_{\omega^k}(t)}{t^2}dt = 0.\tag{2.56}
\end{equation}
In view of (2.56), we may rewrite (2.55) as 
\begin{equation}
     m_{\omega^k}(x) =
\frac{M_{\omega^k}(x)}{x}-\int_x^{\infty}\frac{M_{\omega^k}(t)}{t^2}dt. \tag{2.57}
\end{equation}
From \textit{Theorem 2.1-k}, we see that
\begin{equation}
    \frac{M_{\omega^k}(x)}{x} = \frac{Q_{2,k}^*(\log\log x)}{\log^2x}+\sum_{j=3}^{\nu}\frac{Q_{j,k}^*(\log\log x)}{\log^jx}+O_{\nu}\left(\frac{(\log\log x)^{k\nu+2k}}{\log^{\nu+1}x}\right). \tag{2.58}
\end{equation}
As for the integral on the right in (2.57), if we substitute the expression in \textit{Theorem 2.1-k} for $M_{\omega^k}(t)$ in the integral with $\nu$ replaced $\nu+1$, we will get
\begin{equation}
    -\int_x^{\infty}\frac{M_{\omega^k}(t)}{t^2}dt = -\int_{x}^{\infty}\frac{Q_{2,k}^*(\log\log t)}{t\log^2t}dt+\sum_{j=3}^{\nu+1}\int_x^{\infty}\frac{Q_{j,k}^*(\log\log t)}{t\log^jt}dt+O_{\nu}\left(\int_x^{\infty}\frac{(\log\log t)^{k\nu+3k}}{t\log^{\nu+2}t}dt\right). \tag{2.59}
\end{equation}
In (2.59), the first integral on the right is what gives the leading term in the series expansion of $m_{\omega^k}(x)$. Since each $Q_{j,k}^*(\log\log x)$ is a polynomial in $\log\log x$ of degree $\leq k-1$, with $Q_{2,k}^*(\log\log x)$ of degree $k-1$, the integral expressions on the right of (2.59) are linear combinations of integrals of the form 
\begin{equation}
    I_{r,j}(x):= \int_x^{\infty}\frac{(\log\log t)^r}{t\log^j}dt,\quad \text{where $0 \leq r \leq k-1$, and $j\geq 2$}. \tag{2.60}
\end{equation}
The substitution $\log\log t =v$ converts $I_{r,j}$ into
\begin{equation}
    I_{r,j}(x) = \int_{\log\log x}^{\infty} v^re^{-(j-1)v}dv ,\tag{2.61}
\end{equation}
which can be evaluated by integrating by parts $r$ times to get 
\begin{equation}
    \frac{G_{j,r}(\log\log x)}{\log^{j-1}x}, \tag{2.62}
\end{equation}
where $G_{j,r}$ are polynomials such that
\begin{eqnarray}
    G_{j,r}(x) = \sum_{i=0}^{r} \frac{x^i}{(j-1)^{r-i+1}}. \nonumber
\end{eqnarray}
Combining all such $I_{r,j}$'s arising from (2.59), we get the expression
\begin{equation}
    \frac{F_{1,k}^*(\log\log x)}{\log x} +\frac{F_{2,k}^*(\log\log x)}{\log^2 x}+\cdots  +\frac{F_{\nu,k}^*(\log\log x)}{\log^{\nu} x}\tag{2.63}
\end{equation}
using (2.60), (2,61) and (2.62) if we have $j \leq \nu+1$. Finally, we have
\begin{equation}
    O_{\nu}\left(\int_x^{\infty}\frac{(\log\log t)^{k\nu+3k}}{t\log^{\nu+2}t}dt\right) \ll_{\nu} \frac{(\log\log x)^{k\nu+3k}x}{\log^{\nu+1}x}  .\tag{2.64}
\end{equation}
\textit{Theorem 2.4-k (ii)} follows from (2.57)-(2.64).\qed
\section{The Third Largest Prime Factor}
In what follows, we will provide details of the proofs in the case $k=3$, and sketch how the same methods yield results for all $k>3$.  Thus, we start with the discussion of $P_3(n)$, the third largest prime factor of $n$, but for this, we need results on $P_1(n)$, the largest prime factor of $n$, and $P_2(n)$, the second largest prime factor of $n$, and we begin by stating these. 
\subsection{The counting function for $P_1(n)$}
The fundamental counting function associated with the largest prime factor $P_1(n)$ is defined by 
\begin{eqnarray}
    \Psi_1(x,y) := \Psi(x,y)= \sum_{\substack{n\leq x \\ P_1(n) \leq y}} 1 .
\end{eqnarray}
We let $\alpha = \frac{\log x}{\log y}$. de Bruijn \cite{dB6} has proved that with some $c>0$
\begin{eqnarray}
    \Psi(x,y) \ll xe^{-c\alpha}
\end{eqnarray}
uniformly for $2 \leq y \leq x$. Tenenbaum \cite{Tbook} has shown that $c=\frac{1}{2}$ is admissible in (3.2). de Bruijn \cite{dB5} also proved the following improved bound
\begin{eqnarray}
    \Psi(x,y) \ll x(\log^2 y)e^{-\alpha\log \alpha -\alpha\log \log \alpha +O(\alpha)},
\end{eqnarray}
for $y>\log^2 x$, as well as the asymptotic estimate
\begin{eqnarray}
    \Psi(x,y) \sim x\rho(\alpha) ,
\end{eqnarray}
when $e^{(\log x)^{\frac{3}{5}}} \leq y \leq x$, where $\rho$ satisfies the integro-differential equation
\begin{eqnarray}
    \rho(\alpha) = 1- \int_1^{\alpha} \frac{\rho(u-1)}{u}du ,
\end{eqnarray}
and 
\begin{eqnarray}
    \rho(\alpha) = e^{-\alpha\log \alpha -\alpha\log \log \alpha +O(\alpha)} .
\end{eqnarray}
This shows that $\Psi(x,y)$ is small in comparison to $x$ when $\alpha$ is large. \\ \\
\textbf{Remark:} Even though (3.4) has been improved substantially in terms of the range of values of $y$ (see Hildebrand-Tenenbaum \cite{HT86}), the estimates given in (3.3) and (3.4) suffice for our
purpose. 

\subsection{Ambiguity in the Second and Third Largest Prime Factor}
%Previously, it has been noted that $P_3(n)$ denotes the third largest prime factor of an integer $n$. But how do we really define $P_3(n)$?

The second largest prime factor $P_2(n)$ can be defined in two ways, either as $P_1(n/P_1(n))$, or as the largest prime factor of $n$ strictly less than $P_1(n)$, when $\omega(n)>1$. It was shown in Alladi-Johnson \cite{AJ24}, that there is little difference asymptotically between these two definitions because the number $N_1(x)$ of integers up to $x$ for which $P_1(n)$ repeats is very small in comparison with $x$. More precisely, using (3.3), the following was shown in Alladi-Johnson \cite{AJ24}: \\ \\
\textbf{Lemma 3.1:} \textit{The number $N_1(x)$ of integers $n \leq x$ for which the largest prime factor repeats, satisfies}
\begin{align}
    N_1(x)\ll\frac{x}{e^{\left(\frac{1}{\sqrt 2}+o(1)\right)\sqrt{\log\,x\log\log\,x}}}.
\end{align} 
Thus in Alladi-Johnson \cite{AJ24}, $P_2(n)$ was defined to be the largest prime factor of $n$ that is strictly less than $P_1(n)$ if $\omega(n)\ge 2$, and $P_2(n)=1$, if $\omega(n)<2$. \\ \\
Analogously, one could define $P_3(n)$, the third largest prime factor of $n$, in two ways:
\begin{itemize}
    \item[(i)] using weak inequalities, define $P_2(n)=P_1(n/P_1(n)),$ and $P_3(n)=P_1(n/P_1(n)P_2(n))$, when $\Omega(n)\geq 3$, and $P_3(n)=1$, when $\Omega(n)<3$, or
    \item[(ii)] using strict inequalities, define $P_2(n)$ as the largest prime factor of $n$ which is strictly less than $P_1(n)$, and $P_3(n)$ as the largest prime factor of $n$ which is strictly less than $P_2(n)$, when $\omega(n) \geq 3$, and $P_3(n)=1$, when $\omega(n)<3$.
\end{itemize}
%Of course, there is no trouble for integers with at most 2 prime factors, and for them, we define that $P_3(n)=1$. For the rest, it is easy to understand that we have two different definitions for the same quantity, as follows: 
%\begin{itemize}
%    \item[(i)] $P_3(n):=P_1\left(\frac{n}{P_1(n)P_2(n)}\right)$, where $P_i(n)$ is the $i^{th}$ largest prime factor of $n$. Plainly speaking, we divide the integer by its largest and second largest prime factors, and then define the third largest prime factor of $n$ to be the largest prime factor of the derived quotient. Note here that the second largest prime factor $P_2(n)$ of $n$ is the largest prime factor of $n$ less than $P_1(n)$
 %   \item[(ii)] We can also define the third largest prime factor of $n$ to be the largest prime factor of $n$ less than $P_2(n)$
%\end{itemize}
%Both these definitions are valid, yet there is subtle ambiguity in which one to choose in our context, and to prove our results. At this moment, we observe that the two definitions are consistent with each other, if and only if the integer $n$ does not have repeating largest and second largest prime factors, i.e. the indices of $P_1(n)$ and $P_2(n)$ must be 1 in the prime factorization of $n$. The following lemma will ensure for us that we can work with any of the above definitions as we only need to consider the integers where both of them are consistent, for proving our results. \\ \\
\noindent Analogous to \textit{Lemma 3.1}, we now show that there is little difference asymptotically in the two definitions of $P_3(n)$: \\ \\
\textbf{Lemma 3.2:} \textit{Let $N(x)$ denote the integers less than or equal to $x$ such that either the largest or the second largest prime factors of the integers repeat in their respective prime factorizations. Then}
\begin{eqnarray}
    N(x) \ll \frac{x\log\log x}{\log x}. \nonumber
\end{eqnarray}
\textbf{Proof:} Let us partition the set of integers less than or equal to $x$ with repeating largest or second largest primes into the following sets $N_1(x)$ and $N_2(x)$: we define them as,
\begin{eqnarray}
    N_1(x):=\{n \leq x:\;P_1(n)\; \text{repeats}\} \;\;\;\;\; N_2(x):= \{n \leq x: \;P_1(n) \;\text{does not repeat but $P_2(n)$ does}\} .\nonumber 
\end{eqnarray}
Here, note that there might be integers in $N(x)$ such that $P_1(n)$ and $P_2(n)$ both repeat. By the above definitions, it is automatic that such integers will lie in the set $N_1(x)$. Of course, $N(x)=N_1(x)\cup N_2(x)$. By \textit{Lemma 3.1} above, we already have that
\begin{eqnarray}
N_1(x)\ll\frac{x}{e^{\left(\frac{1}{\sqrt 2}+o(1)\right)\sqrt{\log\,x\log\log\,x}}}.
\end{eqnarray}
Now, for any arbitrary $n \in N_2(x)$, if we denote $P_1(n)=q,P_2(n)=p$, then we have that $n$ is of the order of magnitude $O\left(\frac{x}{p^2}\right)$. Further, by Theorem 6 of \cite{AJ24}, we have that
\begin{eqnarray}
    \Psi_2(x,T) = \sum_{\substack{n \leq x \\ P_2(n) \leq T}} 1 \ll \frac{x\log T}{\log x}.
\end{eqnarray}
Therefore, using (3.8) and (3.9), we have
\begin{eqnarray}
    N_2(x) \ll \Psi_2(x,T) + \sum_{p > T} \frac{x}{p^2} \ll \frac{x\log T}{\log x} + \frac{x}{T\log^2 T}. 
\end{eqnarray}
Now, choosing 
\begin{align}
    T = \frac{\log x}{(\log\log x)^3}, \nonumber
\end{align}
we get our required result. \qed \\ \\
\textbf{Remark:} \textit{Lemma 3.2} shows that the set of integers where the two definitions were inconsistent, is of the order $o(x)$. Therefore, they will be only a part of the error and have no contribution in the main terms of the theorems that follow. So we adopt the definition of $P_3(n)$ using (ii) above, that is using strict inequalities. However, it is to be noted that the bound for $N_1(x)$ in \textit{Lemma 3.1} is much smaller compared to the bound for $N(x)$ in \textit{Lemma 3.2}, and thus the discussion of asymptotics in this paper is more involved compared to that in Alladi-Johnson \cite{AJ24}. \\ \\  
The following lemma is very crucial for proving the uniform distribution of the sequence of the third largest prime factors in reduced residue classes. The lemma will imply that $P_3(n)$ is large for almost all integers, and will also provide a quantitative measure. \\ \\
\textbf{Lemma 3.3} (Tenenbaum \cite{Tbill}) With $\alpha=\frac{\log x}{\log T}$, we have uniformly for $2 \leq T \leq x^{\frac{1}{3}}$
\begin{align}
    \Psi_3(x,T) := \sum_{\substack{n \leq x \\ P_3(n) \leq T}} 1 \ll \frac{x\log\alpha\log T}{\log x} .
\end{align}
\textbf{Remark:} Tenenbaum \cite{Tbill} actually gets a sharper result than what is stated in \textit{Lemma 3.3} and his proof uses intricate analysis. For our purpose, a special case of \textit{Lemma 3.3} namely that
\begin{align}
    \Psi_3(x,T) \ll_{\delta} \frac{x\log\alpha \log T}{\log x}, \quad \text{for $T \leq e^{(\log x)^{1-\delta}}$}
\end{align}
suffices, and this can be established elementarily using the bound for $\Psi(x,y)$ in (3.2). The elementary derivation of a similar bound, more generally, for $\Psi_k(x,T)$, for $k \geq 3$, is given in \S5 (\textit{Theorem 5.2}).

\subsection{Uniform Distribution of the Third Largest Prime Factor of $n$ modulo $\ell$}
\textbf{Theorem 3.4:} \textit{For each integer $\ell \geq 2$, the sequence $\{P_3(n)\}$ of the third largest prime factors of $n$, is uniformly distributed in the reduced residue classes modulo $\ell$. More precisely, for each fixed $\ell \geq 2$, and any $1 \leq j \leq \ell$ satisfying $(j,\ell)=1$, we have that}
\begin{eqnarray}
   N_3(x;\ell,j):=\sum_{\substack{n\leq x\\ P_3(n) \equiv j\;(mod\;\ell)}} 1 = \frac{x}{\varphi(\ell)}+O\left(\frac{x(\log\log x)^3}{\log x}\right) .
\end{eqnarray}
\textbf{Proof:} We note that $P_3(n)=1$ when $\omega(n) \leq 2$. Now,
\begin{eqnarray}
    \sum_{\substack{n \leq x\\ \omega(n)\leq 2}} 1  = \frac{x\log \log x}{\log x} + O\left(\frac{x}{\log x}\right) .
\end{eqnarray}
Therefore, we can just consider the integers with $\omega(n) \geq 3$ in the summation. For a fixed prime $p$, let $S_3(x,p)$ denote the set of integers less than or equal to $x$ with at least 3 distinct prime factors, with square-free first and second largest prime factors and $p$ as the third largest prime factor. Therefore, we have, using \textit{Lemma 3.1} and (3.14)
\begin{eqnarray}
    \sum_{p \leq x^{\frac{1}{3}}} |S_3(x,p)| = [x] - \sum_{\substack{n \leq x\\ \omega(n)\leq 2}} 1= x+O\left(\frac{x\log\log x }{\log x}\right) .
\end{eqnarray}
%Note that (3.28) is possible since we will choose $T$ in such a way that $T$ also goes towards infinity with $x$. \\
Let $N \in S_3(x,p)$. Then, we may write
\begin{align}
    N=m.pqr \quad \text{where} \quad p<q<r\quad \text{are primes and }P_1(m)\leq p.\nonumber
\end{align}
Thus it is clear that $P_1(n)=r$, $P_2(n)=q$ and $P_3(n)=p$. %Recall here that we are only considering integers such that their first and second largest prime number do not repeat. 
Then, with (3.1), we can rewrite $|S_3(x,p)|$ as the following sum:
\begin{eqnarray}
    |S_3(x,p)|= \sum_{\substack{m< \frac{x}{p^3}\\ P_1(m)\leq p}}\; \sum_{p<q< \sqrt{\frac{x}{mp}}}\;\; \sum_{q<r\leq \frac{x}{mpq}} 1  
    = \sum_{p <r \leq \frac{x}{p^2}}\;\; \sum_{p<q<r}\;\; \sum_{\substack{m < \frac{x}{pqr}\\ P_1(m)\leq p}} 1 = \sum_{p <r \leq \frac{x}{p^2}}\;\; \sum_{p<q<r} \Psi\left(\frac{x}{pqr};p\right) .
\end{eqnarray}
Therefore, from (3.15) and (3.16), we get
\begin{eqnarray}
    \sum_{p \leq x^{\frac{1}{3}}} |S_3(x,p)| = \sum_{p \leq x^{\frac{1}{3}}} \;\;\sum_{p <r \leq \frac{x}{p^2}}\;\; \sum_{p<q<r} \Psi\left(\frac{x}{pqr};p\right) = x +O\left(\frac{x\log\log x}{\log x}\right).
\end{eqnarray}
\textbf{Note:} We need to rewrite the triple sum in (3.17) with the innermost sum as a sum over $p$ because we later want to consider the sum over $p$ with the restriction $p \equiv j \;(mod\;\ell)$. Thus, (3.17) is rewritten as 
\begin{align}
    \sum_{p \leq x^{\frac{1}{3}}} |S_3(x,p)| = \sum_{r\leq \frac{x}{6}}\sum_{q<r}^*\sum_{p<q}^*\Psi\left(\frac{x}{pqr};p\right),
\end{align}
where * over the summations means the conditions governing $p,q,r$ are 
\begin{align}
    p\leq x^{\frac{1}{3}},\;p<q<r\leq \frac{x}{6}\quad  \text{and}\quad pqr\leq x. \nonumber
\end{align}
Now, by \textit{Lemma 3.3}, we have for $T \leq e^{\sqrt{\log x}}$,
\begin{eqnarray}
    \sum_{p \leq T} |S_3(x,p)| = O\left(\frac{x\log\alpha\log T}{\log x}\right) .
\end{eqnarray}
Therefore,
(3.17) and (3.19) together give us
\begin{align}
    \Sigma^T :=\sum_{T<p\leq x^{\frac{1}{3}}}|S_3(x,p)| &= \sum_{T<r\leq \frac{x}{T^2}}\sum_{T<q<r}^*\sum_{T<p<q}^*\Psi\left(\frac{x}{pqr};p\right) \nonumber \\
    &= x + O\left( \frac{x\log\log x}{\log x}\right)+ O\left( \frac{x\log \alpha \log T}{\log x}\right).
\end{align}
We now split $\Sigma^T$ as 
\begin{eqnarray}
    \Sigma^T = \Sigma_1^T+\Sigma_2^T+\Sigma_3^T ,
\end{eqnarray}
where
\begin{align}
     \Sigma_1^T &: = \sum_{T<r \leq x^{\frac{1}{3}}} \;\;\sum_{T<q <r}\;\;\sum_{T<p<q} \Psi\left(\frac{x}{pqr},p\right), \hspace{0.5cm} \\
    \Sigma_2^T&:= \sum_{ x^{\frac{1}{3}} < r \leq \frac{x}{T^2}} \;\;\sum_{T<q \leq x^{\frac{1}{3}}}\;\;\sum_{T<p<\min(q,x/qr)}  \Psi\left(\frac{x}{pqr},p\right),  \\
    \Sigma_3^T&:= \sum_{ x^{\frac{1}{3}} < r \leq \frac{x}{T^2}} \;\;\sum_{x^{\frac{1}{3}}<q < r}\;\;\sum_{T<p\leq \frac{x}{qr}}  \Psi\left(\frac{x}{pqr},p\right). 
\end{align}
We now wish to compare $\Sigma^T$ with an expression where the inner sum is replaced by an integral. Replace the inner sums of the above triple sums in (3.22)-(3.24) with corresponding integrals as follows: let us define pseudo-integrals $\mathcal{I}_1^T,\mathcal{I}_2^T$ and $\mathcal{I}_3^T$ as 
\begin{align}
    \mathcal{I}_1^T &: = \sum_{T<r \leq x^{\frac{1}{3}}} \;\;\sum_{T<q <r}\;\;\int_T^q \Psi\left(\frac{x}{pqr},p\right) \frac{dt}{\log t} , \nonumber\\
    \mathcal{I}_2^T&:= \sum_{ x^{\frac{1}{3}} < r \leq \frac{x}{T^2}} \;\;\sum_{T<q \leq x^{\frac{1}{3}}}\;\;\int_T^{\min(q,x/qr)}  \Psi\left(\frac{x}{pqr},p\right)\frac{dt}{\log t},  \nonumber \\
    \mathcal{I}_3^T&:= \sum_{ x^{\frac{1}{3}} < r \leq \frac{x}{T^2}} \;\;\sum_{x^{\frac{1}{3}}<q < r}\;\;\int_T^{\frac{x}{qr}}  \Psi\left(\frac{x}{pqr},p\right)\frac{dt}{\log t} .\nonumber
\end{align}
For $i \in \{1,2,3\}$, we define $E_i = \Sigma_i^T - \mathcal{I}_i^T$. Therefore, from (3.21), we have that
\begin{eqnarray}
    \Sigma^T = \mathcal{I}_1^T+\mathcal{I}_2^T + \mathcal{I}_3^T + E_1+E_2+E_3 .
\end{eqnarray}
We now move ahead to estimate each of the error terms $E_i$. 
\begin{align}
    E_1 &= \sum_{T<r \leq x^{\frac{1}{3}}} \;\;\sum_{T<q <r}\left( \sum_{T<p<q} \Psi\left(\frac{x}{pqr},p\right)-\int_T^q \Psi\left(\frac{x}{pqr},p\right) \frac{dt}{\log t} \right)  \nonumber \\
    &=    \sum_{T<r \leq x^{\frac{1}{3}}} \;\;\sum_{T<q <r}\left( \sum_{T<p<q} \sum_{\substack{n \leq \frac{x}{pqr}\\ P_1(n)\leq p}} 1-\int_T^q \left[\sum_{\substack{n \leq \frac{x}{tqr}\\ P_1(n)\leq p}} 1\right] \frac{dt}{\log t} \right) \hspace{2.4cm} \nonumber \\
    &= \sum_{T<r \leq x^{\frac{1}{3}}} \;\;\sum_{T<q <r}\;\;\sum_{n \leq \frac{x}{Tqr}}\left( \sum_{\max(T,P_1(n))\leq p \leq \min(\frac{x}{nqr},q)} 1 - \int_{\max(T,P_1(n))}^{\min(\frac{x}{nqr},q)} \frac{dt}{\log t}             \right)  .
\end{align}
The difference between the sum and the integral inside the parenthesis on the right of (3.26) can be bounded using the strong form of the Prime Number Theorem. It is important to observe that the error in the strong form of the PNT that we use, namely $\frac{x}{e^{\sqrt{\log x}}}$, is an increasing function of $x$ and therefore, we bound $E_1$ by simply choosing $\frac{x}{nqr}$ to be the upper limit of both the sum and integral and by discarding the lower limit. Thus, we get 
\begin{align}
    |E_1| &\ll \sum_{T<r \leq x^{\frac{1}{3}}} \;\;\sum_{T<q <r}\;\;\sum_{n \leq \frac{x}{Tqr}} \frac{x}{nqr\exp\{\sqrt{\log(\frac{x}{nqr})}\}} \ll \frac{x}{\exp\{\sqrt{\log T}\}}\sum_{T<r \leq x^{\frac{1}{3}}} \;\;\sum_{T<q <r} \frac{\log x}{qr} 
    \nonumber \\
    &\ll \frac{x\log x}{\exp\{\sqrt{\log T}\}}\sum_{T<r \leq x^{\frac{1}{3}}} \frac{\log \log x}{r} \ll \frac{x\log x(\log \log x)^2}{\exp\{\sqrt{\log T}\}}. 
\end{align}
This gives us an estimate for $E_1$. Since the inner sum and the corresponding inner integral in the error $E_2$ are identical to that of $E_1$, we can write $E_2$ as
\begin{eqnarray}
    E_2 = \sum_{ x^{\frac{1}{3}} < r \leq \frac{x}{6}} \;\;\sum_{T<q \leq x^{\frac{1}{3}}}\;\; \sum_{n \leq \frac{x}{Tqr}} \left( \sum_{\max(T,P_1(n))\leq p \leq \min(\frac{x}{nqr},q)} 1 - \int_{\max(T,P_1(n))}^{\min(\frac{x}{nqr},q)} \frac{dt}{\log t}             \right)  .
\end{eqnarray}
We note here that the upper limit in both the sum and integral over $p$ is essentially given by $\min\left(\frac{x}{nqr},q,\frac{x}{qr}\right)$. Since, $n\geq 1$, therefore it suffices to drop $\frac{x}{qr}$. A similar use of the strong form of the Prime Number Theorem therefore yields:
\begin{align}
    |E_2| &\ll \sum_{ x^{\frac{1}{3}} < r \leq \frac{x}{6}} \;\;\sum_{T<q \leq x^{\frac{1}{3}}}\;\; \sum_{n \leq \frac{x}{Tqr}} \frac{x}{nqr\exp\{\sqrt{\log(\frac{x}{nqr})}\}} \ll\frac{x}{\exp\{\sqrt{\log T}\}} \sum_{ x^{\frac{1}{3}} < r \leq \frac{x}{6}} \;\;\sum_{T<q \leq x^{\frac{1}{3}}} \frac{\log x}{qr} \nonumber \\
    &\ll \frac{x\log x}{\exp\{\sqrt{\log T}\}}\sum_{ x^{\frac{1}{3}} < r \leq \frac{x}{6}} \frac{\log\log x}{r} \ll \frac{x\log x(\log\log x)^2}{\exp\{\sqrt{\log T}\}}. 
\end{align}
A similar treatment of $E_3$ therefore will lead us to:
\begin{eqnarray}
    E_3 = \sum_{ x^{\frac{1}{3}} < r \leq \frac{x}{6}} \;\;\sum_{x^{\frac{1}{3}}<q < r} \;\; \sum_{n \leq \frac{x}{Tqr}} \left(        \sum_{\max(T,P_1(n))\leq p \leq \frac{x}{nqr}} 1 - \int_{\max(T,P_1(n))}^{\frac{x}{nqr}} \frac{dt}{\log t}         \right).
\end{eqnarray}
Therefore, again by the strong form of the Prime Number Theorem, we get
\begin{align}
    |E_3| &\ll \sum_{ x^{\frac{1}{3}} < r \leq \frac{x}{6}} \;\;\sum_{x^{\frac{1}{3}}<q < r} \;\; \sum_{n \leq \frac{x}{Tqr}} \frac{x}{nqr\exp\{\sqrt{\log (\frac{x}{nqr})}\}} \ll \frac{x}{\exp\{\sqrt{\log T}\}} \sum_{ x^{\frac{1}{3}} < r \leq \frac{x}{6}} \;\;\sum_{x^{\frac{1}{3}}<q < r} \frac{\log x}{qr} \nonumber \\
    &\ll \frac{x \log x}{\exp\{\sqrt{\log T}\}} \sum_{ x^{\frac{1}{3}} < r \leq \frac{x}{6}} \frac{\log \log x}{r} \ll \frac{x\log x(\log \log x)^2}{\exp\{\sqrt{\log T}\}}. 
\end{align}
Therefore combining all the identities and error estimates from (3.15) to (3.31) , we get that
\begin{eqnarray}
    x = \mathcal{I}_1^T+\mathcal{I}_2^T + \mathcal{I}_3^T + O\left(\frac{x\log\alpha\log T}{\log x}\right) + O\left(\frac{x\log x(\log\log x)^2}{\exp\{\sqrt{\log T}\}}\right).
\end{eqnarray}
Choosing $T = \exp\{(2\log\log x)^2\}$, we finally get that
\begin{eqnarray}
     \mathcal{I}_1^T+\mathcal{I}_2^T + \mathcal{I}_3^T =x+ O\left(\frac{x(\log\log x)^3}{\log x}\right)  .
\end{eqnarray}
We now fix an arbitrary integer $\ell \geq 2$. For any integer $j$ satisfying $ 1 \leq j \leq \ell$ and $(j,\ell)=1$, we consider the set $S_{3,\ell,j}(x)$ of integers $n \leq x$ with non-repeating largest and second largest prime factors with $\omega(n)\geq 2$ and $P_3(n) \equiv j\;(mod\;\ell)$. Then, we have, again from \textit{Lemma 3.2} that
\begin{eqnarray}
    N_3(x;\ell,j) = |S_{3,\ell,j}(x)| + O\left(\frac{x\log\log x}{\log x}\right) .
\end{eqnarray}
We note here that
\begin{eqnarray}
    |S_{3,\ell,j}(x)| = \sum_{\substack{p \leq x^{\frac{1}{3}}\\p\equiv j\;(mod\;\ell)}} |S_3(x,p)| .
\end{eqnarray}
Of course, again using \textit{Lemma 3.3}, we get that
\begin{eqnarray}
    \sum_{\substack{p\leq T\\ p \equiv j\;(mod\;\ell)}} |S_3(x,p)| \leq \sum_{p \leq T} |S_3(x,p)| = O\left(\frac{x\log\alpha\log T}{\log x}\right) .
\end{eqnarray}
We will specify the choice of $T$ in the end of the proof. Therefore, from (3.35) and (3.36), we have
\begin{eqnarray}
    |S_{3,\ell,j}(x)| = \sum_{\substack{
    T<p \leq x^{\frac{1}{3}}\\ p \equiv j\;(mod\;\ell)}} |S_3(x,p)| + O\left(\frac{x\log\alpha \log T}{\log x}\right)+O\left(\frac{x\log\log x}{\log x}\right) .
\end{eqnarray}
We have already set up a summation representation of $|S_3(x,p)|$ in (3.17). Therefore analogously, as done in (3.21), we can write that 
\begin{eqnarray}
     \sum_{\substack{T<p \leq x^{\frac{1}{3}}\\ p \equiv j\;(mod\;\ell)}} |S_3(x,p)| = \Sigma_{1,\ell,j}^T+\Sigma_{2,l\ell,j}^T+\Sigma_{3,\ell,j}^T ,
\end{eqnarray}
where 
\begin{align}
      \Sigma_{1,\ell,j}^T &: = \sum_{T<r \leq x^{\frac{1}{3}}} \;\;\sum_{T<q <r}\;\;\sum_{T<p<q;\;p\equiv j\;(mod\;\ell)} \Psi\left(\frac{x}{pqr},p\right) , \nonumber\\
    \Sigma_{2,\ell,j}^T&:= \sum_{ x^{\frac{1}{3}} < r \leq \frac{x}{6}} \;\;\sum_{T<q \leq x^{\frac{1}{3}}}\;\;\sum_{\substack{T<p<\min(q,x/qr)\\p\equiv j\;(mod\;\ell)}}  \Psi\left(\frac{x}{pqr},p\right), \nonumber \\
    \Sigma_{3,\ell,j}^T&:= \sum_{ x^{\frac{1}{3}} < r \leq \frac{x}{6}} \;\;\sum_{x^{\frac{1}{3}}<q<r}\;\;\sum_{T<p\leq \frac{x}{qr};\;p\equiv j\;(mod\;\ell)}  \Psi\left(\frac{x}{pqr},p\right). \nonumber
\end{align}
As done previously, we will replace these triple summations with corresponding integrals, with the difference being that the integrals are multiplied by a constant $\frac{1}{\varphi(\ell)}$ owing to the PNTAP. Thus, we define $\mathcal{I}_{1,\ell,j}^T,\mathcal{I}_{2,\ell,j}^T$ and $\mathcal{I}_{3,\ell,j}^T$ as 
\begin{align}
    \mathcal{I}_{1,\ell,j}^T &: = \sum_{T<r \leq x^{\frac{1}{3}}} \;\;\sum_{T<q <r}\;\;\frac{1}{\varphi(\ell)}\int_T^q \Psi\left(\frac{x}{pqr},p\right) \frac{dt}{\log t} , \nonumber\\
    \mathcal{I}_{2,\ell,j}^T&:= \sum_{ x^{\frac{1}{3}} < r \leq \frac{x}{6}} \;\;\sum_{T<q \leq x^{\frac{1}{3}}}\;\;\frac{1}{\varphi(\ell)}\int_T^{\min(q,x/qr)}  \Psi\left(\frac{x}{pqr},p\right)\frac{dt}{\log t} , \nonumber \\
    \mathcal{I}_{3,\ell,j}^T&:= \sum_{ x^{\frac{1}{3}} < r \leq \frac{x}{6}} \;\;\sum_{x^{\frac{1}{3}}<q<r}\;\;\frac{1}{\varphi(\ell)}\int_T^{\frac{x}{qr}}  \Psi\left(\frac{x}{pqr},p\right)\frac{dt}{\log t} .\nonumber
\end{align}
We aim to replace $\Sigma_{i,\ell,j}^T$ with $\mathcal{I}_{i,\ell,j}^T$, $i \in \{1,2,3\}$, in (3.38). So we denote the difference by
\begin{align}
    E_{i,\ell,j}= \Sigma_{i,\ell,j}^T - \mathcal{I}_{i,\ell,j}^T,\quad i=1,2,3. \nonumber
\end{align}
%Our current case is not very different from the one we have already taken care of before, with just a modulo condition on the prime $p$ and a constant factor. So, following similar calculations as done above, we can get the following expressions for each of the error terms:
Thus,
\begin{align}
    E_{1,\ell,j} &= \sum_{T<r \leq x^{\frac{1}{3}}} \;\;\sum_{T<q <r}\;\;\sum_{n \leq \frac{x}{Tqr}}\left( \sum_{\substack{\max(T,P_1(n))\leq p \leq \min(\frac{x}{nqr},q)\\ p\equiv j\;(mod\;\ell)}} 1 - \frac{1}{\varphi(\ell)}\int_{\max(T,P_1(n))}^{\min(\frac{x}{nqr},q)} \frac{dt}{\log t}             \right), \nonumber \\
     E_{2,\ell,j} &= \sum_{ x^{\frac{1}{3}} < r \leq \frac{x}{6}} \;\;\sum_{T<q \leq x^{\frac{1}{3}}}\;\; \sum_{n \leq \frac{x}{Tqr}} \left( \sum_{\substack{\max(T,P_1(n))\leq p \leq \min(\frac{x}{nqr},q)\\p\equiv j\;(mod\;\ell)}} 1 - \frac{1}{\varphi(\ell)}\int_{\max(T,P_1(n))}^{\min(\frac{x}{nqr},q)} \frac{dt}{\log t}             \right) , \nonumber  \\
       E_{3,\ell,j}& = \sum_{ x^{\frac{1}{3}} < r \leq \frac{x}{6}} \;\;\sum_{x^{\frac{1}{3}}<q<r} \;\; \sum_{n \leq \frac{x}{Tqr}} \left(        \sum_{\substack{\max(T,P_1(n))\leq p \leq \frac{x}{nqr}\\ p\equiv j\;(mod\;\ell)}} 1 - \frac{1}{\varphi(\ell)}\int_{\max(T,P_1(n))}^{\frac{x}{nqr}} \frac{dt}{\log t}         \right) .\nonumber
\end{align}
We bound $E_{i,\ell,j}$, for $i=1,2,3$, in exactly the same form as we bounded $E_1,E_2,E_3$ above, except that here we use the strong form of the Prime Number Theorem in Arithmetic Progressions, and the fact that the error $\frac{x}{e^{\sqrt{\log x}}}$ is an increasing function of $x$. So, we get that for every $i$,
\begin{eqnarray}
    E_{i,\ell,j} \ll \frac{x\log x(\log\log x)^2}{\exp\{\sqrt{\log T}\}} .\nonumber 
\end{eqnarray}
Using this above estimate in (3.38), we therefore get that
\begin{eqnarray}
     \sum_{\substack{T<p \leq x^{\frac{1}{3}}\\ p \equiv j\;(mod\;\ell)}} |S_3(x,p)| = \mathcal{I}_{1,\ell,j}^T+\mathcal{I}_{2,\ell,j}^T+\mathcal{I}_{3,\ell,j}^T + O\left(\frac{x\log x(\log \log x)^2}{\exp\{\sqrt{\log T}\}}\right).
\end{eqnarray}
But, note that $\mathcal{I}_{i,\ell,j}^T=\frac{1}{\varphi(\ell)} \mathcal{I}_i^T$. Therefore, combining equations (3.34)-(3.39) along with this fact, we get that
\begin{eqnarray}
N_3(x;\ell,j) = \frac{1}{\varphi(\ell)}(\mathcal{I}_{1}^T+\mathcal{I}_{2}^T+\mathcal{I}_{3}^T) +  O\left(\frac{x\log\alpha\log T}{\log x}\right) + O\left(\frac{x\log x(\log\log x)^2}{\exp\{\sqrt{\log T}\}}\right) .
\end{eqnarray}
Taking $T = \exp\{(2\log\log x)^2\}$, we get from (3.40) that 
\begin{eqnarray}
    N_3(x;\ell,j) = \frac{1}{\varphi(\ell)}(\mathcal{I}_{1}^T+\mathcal{I}_{2}^T+\mathcal{I}_{3}^T) + O\left(\frac{x(\log\log x)^3}{\log x}\right).
\end{eqnarray}
Finally, from (3.33) and (3.41), we finally get our desired result:
\begin{eqnarray}
    N_3(x;\ell,j) = \frac{x}{\varphi(\ell)} + O\left(\frac{x(\log\log x)^3}{\log x}\right). \nonumber
\end{eqnarray}\qed \\ \\
\textbf{Remark 3.4.1:} The idea to utilize the identity like (3.17) and the replacement of the inner sum in (3.18) by an integral goes back to \cite{KA1977} and indeed, has been used in all subsequent papers (\cite{Da17}, \cite{SW19}, \cite{Se25} and so on) on applications of this duality identity where such an uniform distribution result is needed. \\ \\
\textbf{Remark 3.4.2:} Here we have used the de la Vallee Poussin error term in the strong form of the Prime Number Theorem, namely 
\begin{align}
    \frac{x}{e^{\sqrt{\log x}}}.\nonumber
\end{align}
If instead, we use stronger error terms such as 
\begin{align}
    \frac{x}{e^{(\log x)^{\delta}}},\quad \text{where $\delta >\frac{1}{2}$}, \nonumber
\end{align}
then the error term in \textit{Theorem 3.4} will improve to 
\begin{align}
    \frac{x(\log\log x)^{1+\frac{1}{\delta}}}{\log x}. \nonumber
\end{align} 
\textit{Theorem 3.4} is of extreme importance in this paper, as it paves the path for our main result that will come in the next section. We will prove three more theorems, building up to our final one. Note here that very interestingly, the density of the integers with their first, second, and third largest prime factors in reduced residue classes $modulo \;\ell$ separately, are the same, i.e. $\frac{1}{\varphi(\ell)}$. Later, we will see that it is also the case for the distribution of $k^{th}$ largest prime factors in reduced residue classes, for every $k$.

\section{Main Results for $k=3$ using Third Order Duality}
We start with the following theorem that give us an estimate of the sum of $\mu(n)\omega(n)^2$ with the condition that the smallest prime factor of $n$ is congruent to $j\;(mod\;\ell)$. We need to use the result in \textit{Theorem 3.4} together with the third order duality to establish the following theorem. \\ \\
\textbf{Theorem 4.1} \textit{For integers $j,\ell$ satisfying $1 \leq  j \leq \ell$ and $(j,\ell)=1$, we have}
\begin{eqnarray}
   M_{\omega^2}(x,j,\ell):= \sum_{\substack{n \leq x\\p_1(n) \equiv j\;(mod\;\ell)}} \mu(n)\omega(n)^2 \ll \frac{x(\log\log x)^5}{\log x}.
\end{eqnarray}
\textbf{Proof:} We use the Duality Lemma with $k=3$. This yields
\begin{align}
      \sum_{1<d|n}\mu(d){\omega(d)-1 \choose 2}f(p_1(d))=-f(P_{3}(n)) .
\end{align}
Next, we apply (4.2) by choosing $f$ to be an arithmetic function on primes defined as follows:
\begin{eqnarray}
    f(n) = \begin{cases}
        1, \;\text{if $n$ is a prime and $n \equiv j\;(mod\;\ell)$} \\
        0,\;\text{otherwise.}
    \end{cases}
\end{eqnarray}
Now, using M\"obius inversion in (4.2), we have that
\begin{eqnarray}
  \frac{1}{2}\sum_{n \leq x} \mu(n)(\omega(n)-1)(\omega(n)-2)f(p_1(n))= -  \sum_{n\leq x} \sum_{d|n} \mu\left(\frac{n}{d}\right) f(P_3(d)).
\end{eqnarray}
Applying the hyperbola method in the double sum of the RHS of (4.4), we get
\begin{eqnarray}
    \sum_{n\leq x} \sum_{d|n} \mu\left(\frac{n}{d}\right) f(P_3(d)) = \sum_{m \leq T}\mu(T)\sum_{d \leq \frac{x}{m}} f(P_3(d))+\sum_{d \leq \frac{x}{T}}f(P_3(d))\sum_{T<m\leq \frac{x}{d}}\mu(m) =: \mathcal{S}_1+\mathcal{S}_2 ,
\end{eqnarray}
where $\mathcal{S}_1$ and $\mathcal{S}_2$ are respectively the double sums in the RHS of (4.5).
Looking at the double sum $\mathcal{S}_2$, we have
\begin{eqnarray}
    \mathcal{S}_2 = \sum_{d \leq \frac{x}{T}}f(P_3(d))\sum_{T<m\leq \frac{x}{d}}\mu(m) \ll \sum_{d \leq \frac{x}{T}}f(P_3(d)) \frac{x}{de^{\sqrt{\log (x/d)}}}\ll \frac{x\log x}{e^{\sqrt{\log T}}} .
\end{eqnarray}
Further, looking at the inner sum of $\mathcal{S}_1$, we observe that using \textit{Theorem 3.4}
\begin{eqnarray}
    \sum_{d \leq \frac{x}{m}} f(P_3(d)) = \frac{x}{m\varphi(\ell)} + O\left(\frac{x(\log\log x)^3}{m\log(x/m)}\right) .\nonumber
\end{eqnarray}
Therefore continuing with the estimating of $\mathcal{S}_1$ using the above, we have
\begin{align}
\mathcal{S}_1&=  \sum_{m \leq T}\mu(m)\sum_{d \leq \frac{x}{m}} f(P_3(d)) = \frac{x}{\varphi(\ell)} \sum_{m \leq T}\frac{\mu(m)}{m}+O\left(x(\log\log x)^3\sum_{m\leq T}\frac{1}{m\log(x/m)}\right)\nonumber \\ 
&\ll \frac{x}{\varphi(\ell)e^{\sqrt{\log (T)}}}+O\left(\frac{x\log T(\log\log x)^3}{\log (x/T)}\right) .
\end{align}
Now, by the following suitable choice of $T=e^{(2\log\log x)^2}$, we get from (4.5)-(4.7) that
\begin{eqnarray}
     \sum_{n\leq x} \sum_{d|n} \mu\left(\frac{n}{d}\right) f(P_3(d)) \ll \frac{x(\log\log x)^5}{\log x} .
\end{eqnarray}
But, here we note that
\begin{align}
    \sum_{n \leq x} \mu(n)(\omega(n)-1)(\omega(n)-2)f(p_1(n)) = \sum_{n \leq x} \mu(n)\omega(n)^2f(p_1(n))-3\sum_{n \leq x} \mu(n)\omega(n)f(p_1(n))\nonumber \\
    +2\sum_{n \leq x} \mu(n)f(p_1(n)), \nonumber
\end{align}
which gives us that
\begin{align}
    M_{\omega^2}(x,j,\ell) =  \sum_{n \leq x} \mu(n)(\omega(n)-1)(\omega(n)-2)f(p_1(n))+3\sum_{n \leq x} \mu(n)\omega(n)f(p_1(n)) \nonumber \\
    -2\sum_{n \leq x} \mu(n)f(p_1(n)).
\end{align}
Also, from \textit{Theorem 8} in \cite{AJ24} and \textit{Theorem 3} in \cite{KA1977}, we have
\begin{eqnarray}
    \sum_{n\leq x}\mu(n)\omega(n)f(p_1(n))\ll \frac{x(\log\log x)^4}{\log x}\;\;\;\;\text{and}\;\;\;\; \sum_{n\leq x}\mu(n)f(p_1(n)) \ll xe^{-(\log x)^{1/3}} .
\end{eqnarray}
Therefore, from (4.4), along with (4.8), (4.9) and (4.10), we finally get our required result. \qed \\ \\
The importance of \textit{Theorem 3.4} has already been realized in the above theorem. The fact that the average of the third largest prime factors of integers in reduced residue class $j\;(mod \;\ell)$ is $\frac{1}{\varphi(\ell)}$ will also prove to be important in one of the upcoming theorems. \\ \\
\textbf{Theorem 4.2:} \textit{For integers $j,\ell$ satisfying $1 \leq j \leq \ell$ and $(j,\ell)=1$, we have}
\begin{eqnarray}
    \sum_{\substack{n \leq x\\ p_1(n)\equiv j\;(mod\;\ell)}}\mu(n)\omega(n)^2\left\{\frac{x}{n}\right\}\ll\frac{x(\log\log x)^{7/2}}{\sqrt{\log x}}. \nonumber
\end{eqnarray}
\textbf{Proof:} We will prove \textit{Theorem 4.2} using \textit{Theorem A} and therefore, a correct choice of the functions $\eta(x)$
 and $\beta(x)$ will suffice. 
We choose the sequence $\{a_n\}_n$ such that for all positive integers $n$, $a_n=\mu(n)\omega(n)^2f(n)$, where $f$ is defined in (4.3). Then, by \textit{Theorem 4.1}, we choose $\eta(x)$ in \textit{Theorem A} to be
\begin{align}
    \eta(x) = \frac{(\log\log x)^5}{\log x} .
\end{align}
Of course, $\eta(x)$ satisfies the conditions in (2.34a) and (2.34b). Further, we have that
\begin{align}
    \sum_{n\leq x}|a_n| \leq \sum_{n\leq x}\omega(n)^2 \ll x(\log\log x)^2 .
\end{align}
Therefore, from (4.12), we choose $\beta(x)$ in \textit{Theorem A} to be
\begin{align}
    \beta(x) = (\log\log x)^2.
\end{align}
Again, $\beta(x)$ satisfies the conditions in (2.36). Thus, by \textit{Theorem A}, we conclude that 
\begin{align}
      \sum_{\substack{n \leq x\\ p_1(n)\equiv j\;(mod\;\ell)}}\mu(n)\omega(n)^2\left\{\frac{x}{n}\right\} =  \sum_{n \leq x}a_n\left\{\frac{x}{n}\right\}\ll\frac{x(\log\log x)^{7/2}}{\sqrt{\log x}},
\end{align}
which proves our theorem. \qed \\ \\
We next consider the sum with $\left[\frac{x}{n}\right]$ replacing $\left\{\frac{x}{n}\right\}$. We then have the following theorem: \\ \\
\textbf{Theorem 4.3} \textit{For integers $j,\ell$ satisfying $1 \leq j \leq \ell$ and $(j,\ell)=1$, we have} 
\begin{eqnarray}
     \sum_{\substack{n \leq x\\ p_1(n)\equiv j\;(mod\;\ell)}}\mu(n)\omega(n)^2\left[\frac{x}{n}\right]\ll \frac{x(\log\log x)^3}{\log x}.
\end{eqnarray}
\textbf{Proof:} We recall the function $f$ defined in (4.3). We replace the condition in the summation using $f$ to get
\begin{eqnarray}
      \sum_{\substack{n \leq x\\ p_1(n)\equiv j\;(mod\;\ell)}}\mu(n)\omega(n)^2\left[\frac{x}{n}\right] = \sum_{n \leq x}\mu(n)\omega(n)^2f(p_1(n))\left[\frac{x}{n}\right] = \sum_{n \leq x} \sum_{d|n}\mu(d)\omega(d)^2f(p_1(d)) .
\end{eqnarray}
Rewriting $\omega(d)^2$ as follows,
\begin{align}
    \omega(d)^2 = (\omega(d)-1)(\omega(d)-2)+3(\omega(d)-1)+1, \nonumber
\end{align}
we get from (4.16) that
\begin{align}
     \sum_{\substack{n \leq x\\ p_1(n)\equiv j\;(mod\;\ell)}}\mu(n)\omega(n)^2\left[\frac{x}{n}\right] = \sum_{n \leq x}& \sum_{d|n}\mu(d)(\omega(d)-1)(\omega(d)-2)f(p_1(d))\left[\frac{x}{n}\right]\nonumber \\
     &+3\sum_{n \leq x} \sum_{d|n}\mu(d)(\omega(d)-1)f(p_1(d))\left[\frac{x}{n}\right]
     +\sum_{n \leq x} \sum_{d|n}\mu(d)f(p_1(d))\left[\frac{x}{n}\right].
\end{align}
Using the Duality Identity in (1.11) for $k=3$, $k=2$ and $k=1$ respectively, for the inner sums in each of the double sums on the right of (4.17), we get that
\begin{eqnarray}
     \sum_{\substack{n \leq x\\ p_1(n)\equiv j\;(mod\;\ell)}}\mu(n)\omega(n)^2\left[\frac{x}{n}\right] = -2\sum_{n\leq x}f(P_3(n))+3\sum_{n\leq x}f(P_2(n)) -\sum_{n\leq x}f(P_1(n)).
\end{eqnarray}
Therefore, using \textit{Theorem 3.4} above, \textit{Theorem 7} in \cite{AJ24} and \textit{Theorem 4} in \cite{KA1977} together in (4.18), we deduce that
\begin{align}
         \sum_{\substack{n \leq x\\ p_1(n)\equiv j\;(mod\;\ell)}}\mu(n)\omega(n)^2\left[\frac{x}{n}\right]
         &= (-2+3-1)\frac{x}{\varphi(\ell)}+ O\left(\frac{x(\log\log x)^3}{\log x}\right)\nonumber \\
         &=O\left(\frac{x(\log\log x)^3}{\log x}\right).
\end{align}
Hence, we have our theorem proved. \qed \\ \\
The beauty of \textit{Theorem 4.3} is that the surprising cancellation of the main term that takes place after the use of the Duality Lemma. This is why, to have exact same averages of all the largest, second largest and third largest prime factors in the same reduced residue classes is so significant in this proof. We now go on to state and prove the final theorem of the section, which is our main result for the case $k=3$. We have: \\ \\
\textbf{Theorem 4.4:} \textit{For integers $j,\ell$ satisfying $1 \leq j \leq \ell$ and $(j,\ell)=1$, we have} 
\begin{eqnarray}
    \sum_{\substack{n \leq x \\ p_1(n)\equiv j\;(mod\;\ell)}} \frac{\mu(n)\omega(n)^2}{n} = O\left(\frac{(\log\log x)^{7/2}}{\sqrt{\log x}}\right).
\end{eqnarray}
\textit{Hence, taking limit $x\to \infty$, we have that}
\begin{eqnarray}
     \sum^{\infty}_{\substack{2 \leq n \\ p_1(n)\equiv j\;(mod\;\ell)}} \frac{\mu(n)\omega(n)^2}{n} =  0.
\end{eqnarray}
\textbf{Proof:} Adding inequations (4.14) and (4.15) of \textit{Theorem 4.2} and \textit{Theorem 4.3}, we get that
\begin{align}
    \sum_{\substack{n \leq x \\ p_1(n)\equiv j\;(mod\;\ell)}}\mu(n)\omega(n)^2\cdot \frac{x}{n} &=  \sum_{\substack{n \leq x \\ p_1(n)\equiv j\;(mod\;\ell)}}\mu(n)\omega(n)^2\left\{\frac{x}{n}\right\}+ \sum_{\substack{n \leq x \\ p_1(n)\equiv j\;(mod\;\ell)}}\mu(n)\omega(n)^2\left[\frac{x}{n}\right] \nonumber \\
    &\ll \frac{x(\log\log x)^{7/2}}{\sqrt{\log x}} + \frac{x(\log\log x)^3}{\log x} \ll \frac{x(\log\log x)^{7/2}}{\sqrt{\log x}} .
\end{align}
Canceling by a factor of $x$ on both sides of (4.22), we therefore have
\begin{eqnarray}
    \sum_{\substack{n \leq x \\ p_1(n)\equiv j\;(mod\;\ell)}} \frac{\mu(n)\omega(n)^2}{n} = O\left(\frac{(\log\log x)^{7/2}}{\sqrt{\log x}}\right). \nonumber
\end{eqnarray}
Finally, taking the limit $x \to \infty$, we get (4.21). \qed \\ \\
We now state our arithmetic density version of \textit{Theorem 4.4} as its corollary. One can follow the same methods used above to get a quantitative evaluation of the same. \\ \\
\textbf{Corollary 4.5:} \textit{For integers $j,\ell$ satisfying $1 \leq j \leq \ell$ and $(j,\ell)=1$, we have that}
\begin{eqnarray}
       \sum^{\infty}_{\substack{2 \leq n \\ p_1(n)\equiv j\;(mod\;\ell)}} \frac{\mu(n){\omega(n)-1 \choose 2}}{n} =  -\frac{1}{\varphi(\ell)}.
\end{eqnarray}
\vspace{0.5cm} \\ \\
\textbf{Remark:} This is only a part, but an important one, of the actual goal of this paper. As said earlier, we aim to prove that
\begin{eqnarray}
    \sum_{\substack{n=2 \\ p_1(n) \equiv j\;(mod\; \ell)}}^{\infty} \frac{\mu(n)\omega(n)^{k-1}}{n} = 0.
\end{eqnarray}
%The $k=2$ case is proved in [A-J 2024] and above, we have proved the $k=3$ case. Note here that since we will prove that (4.24) holds for all integers $k \geq 2$, it does not make any difference in writing $k$ instead of $k-1$, and rather the former seems to be a better option. But we still insist on writing $k-1$ as the exponent of $\omega(n)$, to make the case consistent to the order of duality that we use to prove the theorem. For example, Alladi and Johnson used the second order duality to prove their result, and hence we call it the $k=2$ case. We have used the third order duality above, and hence, we name it the $k=3$ case. Therefore, writing $k-1$ in the exponent of $\omega(n)$ eradicates any confusion regarding the ambiguity of $k$. \\ \\
We now move on to get the general form of the results that we have already proved in the $k=3$ case. It is important to note here that the techniques used to prove results in \cite{KA1977} were drastically different, and the error estimates were much stronger. We will see, as we continue, that as the value of $k$ increases, the error estimates in every step become worse, but interestingly, are still small enough to give us our required theorems. For every $k \geq 3$, we use almost exactly the same technique, and hence, the use of Principle of Mathematical Induction would be our main tool to reach our goal. We would of course show the key steps that are required to prove the general theorems for the reader to get a clear picture, but would mostly refer to the proofs done in the above sections.

\section{The $k^{th}$ largest Prime Factor and its Uniform Distribution in Reduced Residue Classes $mod \; \ell$}
As we saw in the case of the third largest prime factor $P_3(n)$, we can also define the $k^{th}$ largest prime factor, denoted as $P_k(n)$, in two ways:
\begin{itemize}
    \item[(i)] $P_j(n):=P_1\left(\frac{n}{\prod_{i=1}^{j-1}P_i(n)}\right)$, $j=2,\cdots  ,k$, if $\Omega(n)\geq k$, and $P_k(n)=1$, if $\Omega(n)<k$.
    \item[(ii)] If $n=p_1^{e_1}\cdot p_2^{e^2}\cdots   p_r^{e_r}$, with $p_1>p_2>\cdots  >p_r$, then $P_j(n):=p_k$, for $j=1,\cdots   ,k$, if $r= \omega(n)\geq k$, and if $r<k$, then $P_k(n)=1$.
\end{itemize}
It is quite clear that the definitions are different for integers $n$ whose prime factors would repeat. But, as done before, we bypass this inconsistency with the help of the following lemma, which generalizes \textit{Lemma 3.1}. \\ \\
\textbf{Lemma 5.1:} \textit{For an integer $k \geq 3$, let $N_{k-1}(x)$ denote the integers less than or equal to $x$ such  at least one of the $P_i(n)$'s, for $1 \leq i \leq k-1$, repeats in the prime factorization of $n$. Then }
\begin{eqnarray}
    N_{k-1}(x) \ll_k \frac{x(\log\log x)^{k-2}}{\log x}.
\end{eqnarray}
\textbf{Proof:} We start by considering the function that counts the number of integers $\leq x$ whose $k^{th}$ largest prime factor in the sense of (ii) is $\leq T$. 
\begin{align}
    \Psi_{k}(x,T) := \sum_{\substack{n\leq x \\P_k(n)\leq T}}1 .
\end{align}
It follows from a result of Tenenbaum \cite{Tbill} that for
%To prove this lemma, we will use the principle of mathematical induction on a result that we assume here and prove as a theorem (see \textit{Theorem 5.2}) following this lemma. We assume that for some 
$T\leq \exp^{(\log x)^{1-\delta}},$ where $\delta>0$, and for an arbitrary fixed integer $k \geq 2$,
\begin{eqnarray}
  \Psi_k(x,T) \ll_{k,\delta}\frac{x(\log\log x)^{k-2}\log T}{\log x}.
\end{eqnarray}
We observe that \textit{Lemma 3.3} in \S3.2 is the case $k = 3$ (also see \textit{Theorem 6} in \cite{AJ24} for the $k=2$ case). So, we now assume that $k\geq 3$. \\ \\ %Then, for $T \leq e^{(\log x)^{1-\delta}}$, for some $\delta>0$, we have that 
%\begin{align}
 %   \Psi_{k-1}(x,T) \ll \frac{x(\log\log x)^{k-3}\log T}{\log x}.
%\end{align}
%On the other hand, if we choose an arbitrary $n \in N_{k-1}(x)$ such that $P_k(n)=p$, then by our hypothesis, the number of such integers $n$ is of the order of magnitude $O\left((k-1)\frac{x}{p^2}\right)$. Therefore, using (5.4), we have that
To bound the number of integers $\leq x$ for which at least one of the prime factors $P_1(n), P_2(n),\cdots  ,P_{k-1}(n)$ repeats, we consider two cases:  \\ \\
\underline{Case-I:} $P_{k-1}(n)\leq T$. \\ 
In this case, the number of such integers is trivially 
\begin{align}
    \leq \Psi_{k-1}(x,T) \ll_{k,\delta} \frac{x(\log\log x)^{k-3}\log T}{\log x}.
\end{align}
\underline{Case-II:} $P_{k-1}>T$. \\
In this case, each of $P_1(n),P_2(n),\cdots  , P_{k-1}(n)$ is $>T$. Now, the number of integers $\leq x$ for which a given prime factor $p$ repeats is $O\left(\frac{x}{p^2}\right)$. Thus, considering the $k-1$ cases of $P_1(n),P_2(n),\cdots,P_{k-1}(n)$ repeating, we deduce that the number of integers under Case-II, is 
\begin{align}
    \leq \sum_{p>T}(k-1)\cdot \frac{x}{p^2}\ll (k-1)\cdot 
    \frac{x}{T\log^2T}
\end{align}
Thus,
\begin{align}
    N_{k-1}(x)\ll_{k,\delta} \frac{x(\log\log x)^{k-3}\log T}{\log x}+ (k-1)\cdot 
    \frac{x}{T\log^2T}. 
\end{align}
With the choice of 
\begin{align}
    T = \frac{\log x}{(\log\log x)^k}, \nonumber
\end{align}
we have our required bound from (5.6). As the choice of $k$ is arbitrary, our lemma is proved. \qed \\ \\
\textbf{Remark I:} We see that the \textit{Lemma 5.1} makes it clear that the integers where the above definitions of $P_k(n)$ are inconsistent, whose size is small as compared to $x$. Henceforth, the $P_k(n)$ will be defined in a strict sense as in (ii). \\ \\
\textbf{Remark II:} Tenenbaum \cite{Tbill} established a stronger version of (5.3) valid uniformly for $T \leq x^{\frac{1}{k}}$, with $\log\log x$ replaced by $\log\log x-\log\log T$. His proof used analytic techniques. For our purpose, since, $T\leq e^{(\log x)^{1-\delta}}$, we have stated the upperbound in (5.3) with $\log\log x$ in place of $\log\log x-\log\log T$. We now give an elementary proof of (5.3) in the stated range. \\ \\
%We will now prove the assumption that we used to prove \textit{Lemma 5.1}. We will state the result as a theorem, and it will prove to be extremely important to show the uniform distribution of $P_k(n)$ in the reduced residue class $mod\;\ell$. It might be quite understandable now that the following lemma is again an elementary generalization of \textit{Lemma 3.3}, which played an important role in the proof of the uniform distribution of $P_3(n)$. \\ \\
\textbf{Theorem 5.2:} \textit{$P_k(x) \rightarrow \infty$ as $x \rightarrow \infty$ almost always, for any positive integer $k$. Further, for a fixed positive integer $k \geq 2$, we have
\begin{eqnarray}
  \Psi_k(x,T)=  \sum_{\substack{n \leq x\\P_k(n)\leq T}} 1 \ll_{k,\delta} \frac{x (\log\log x)^{k-2} \log T}{\log x} ,
\end{eqnarray}
for all $T \leq e^{(\log x)^{1-\delta}}$, where $\delta >0$}. \\ \\
\textbf{Proof:}  We start with the following routine counting functions:
\begin{eqnarray}
    \Psi_k(x,T) = \sum_{\substack{n \leq x\\ P_k(n) \leq T}} 1 \;\;\;\;\text{and}\;\;\;\; \Psi_k^{(p)}(x) = \sum_{\substack{n \leq x\\P_k(n) = p}} 1   .
\end{eqnarray}
Therefore, we have
\begin{eqnarray}
      \Psi_k(x,T) = \sum_{p \leq T} \Psi_k^{(p)}(x) .
\end{eqnarray}
We first estimate $\Psi_k^{(p)}(x)$. Let $n=m.pp_{k-1}p_{k-2}\cdots  p_2p_1$ and $P_i(n)=p_i$, where $1 \leq i \leq k-1 $, $P_k(n)=p\leq e^{(\log x)^{1-\delta}}$ and $P_1(m) \leq p$. Note here that, $p<p_{k-1}<p_{k-2}<\cdots  <p_2<p_1$. For a fixed $p$, we have
\begin{align}
    \Psi_k^{(p)}(x) &= \sum_{\substack{m \leq \frac{x}{p^k}\\ P_1(m)\leq p}}\;\; \sum_{p^{k-1} \leq p_{k-1}p_{k-2}\cdots  p_2p_1 \leq \frac{x}{mp}} 1 \nonumber\\ 
    &\ll \sum_{\substack{m \leq \frac{x}{p^k}\\ P_1(m)\leq p}} \frac{x \left(\log \log(\frac{x}{mp})\right)^{k-2}}{(k-2)!mp\log(\frac{x}{mp})},\;\;\left[\text{since $\pi_{k-1}(x)\ll\frac{x(\log \log x)^{k-2} }{(k-2)!\log x}$}\right] .
\end{align}
We split the summation on the right of (5.10) as follows:
\begin{eqnarray}
     \Psi_k^{(p)}(x) \ll \sum_{\substack{m \leq \exp{(\log x)^{1-\frac{\delta}{2}}}\\ P_1(m)\leq p}} \frac{x \left(\log \log(\frac{x}{mp})\right)^{k-2}}{mp\log(\frac{x}{mp})} + \sum_{\substack{\exp{(\log x)^{1-\frac{\delta}{2}}} \leq m \leq \frac{x}{p^k}\\P_1(m)\leq p}} \frac{x \left(\log \log(\frac{x}{mp})\right)^{k-2}}{mp\log(\frac{x}{mp})} \;\; =: \; \mathcal{S}_{1} + \mathcal{S}_{2} .
\end{eqnarray}
\textbf{\underline{CASE-I:}} Estimating $\mathcal{S}_{1}$:\\
We observe that $(\log \log(x/mp))^{k-2}\ll(\log \log x)^{k-2} $ and $\log (x/mp) \gg \log x$. Using these, we have from $\mathcal{S}_1$ that 
\begin{align}
    \mathcal{S}_{1} \ll \frac{x(\log \log x)^{k-2} }{p\log x} \sum_{\substack{m \leq \exp{(\log x)^{1-\frac{\delta}{2}}}\\P_1(m)\leq p}} \frac{1}{m} \ll \frac{x(\log \log x)^{k-2} }{p \log x} \prod_{\substack{q \leq p \\ q\;\text{is a prime}}} \left(1-\frac{1}{q}\right)^{-1} \ll \frac{x (\log \log x)^{k-2}  \log p}{p\log x} .
\end{align}
\textbf{\underline{CASE-II:}} Estimating $\mathcal{S}_{2}$: \\
Here, using $\log (x/mp)\geq \log p$, we get
\begin{eqnarray}
    \mathcal{S}_{2} \leq \frac{x(\log \log x)^{k-2} }{p\log p} \sum_{\substack{\exp{(\log x)^{1-\frac{\delta}{2}}} \leq m \\ P_1(m)\leq p}} \frac{1}{m} < \frac{x(\log \log x)^{k-2} }{p\log p} \int_{\exp(\log x)^{1-\frac{\delta}{2}}}^{\infty} \frac{d \Psi(t,p)}{t} .
\end{eqnarray}
Now, looking at the integral above, we use integration by parts to get
\begin{eqnarray}
     \int_{\exp(\log x)^{1-\frac{\delta}{2}}}^{\infty} \frac{d \Psi(t,p)}{t} = \frac{\Psi(t,p)}{t} \Big|_{\exp{(\log x)^{1-\frac{\delta}{2}}}}^{\infty} +  \int_{\exp(\log x)^{1-\frac{\delta}{2}}}^{\infty} \frac{\Psi(t,p)}{t^2}dt.
\end{eqnarray}
Using the estimate (3.2) of $\Psi(x,T)$ in (5.14), we have
\begin{eqnarray}
      \int_{\exp(\log x)^{1-\frac{\delta}{2}}}^{\infty} \frac{d \Psi(t,p)}{t} \ll e^{-\alpha_0/2}+   \int_{\exp(\log x)^{1-\frac{\delta}{2}}}^{\infty} \frac{dt}{te^{\alpha(t)/2}},
\end{eqnarray}
where 
\begin{eqnarray}
    \alpha_0 > \frac{\log e^{(\log x)^{1-\frac{\delta}{2}}} }{\log e^{{(\log x)}^{1-\delta} }} = (\log x)^{\frac{\delta}{2}} \;\;\;\text{and} \;\;\;\alpha(t) = \frac{\log t}{\log p}>\frac{\log e^{(\log x)^{1-\frac{\delta}{2}}} }{\log e^{{(\log x)}^{1-\delta} }} = (\log x)^{\frac{\delta}{2}}. \nonumber 
\end{eqnarray}
This helps us rewrite (5.15) as follows:
\begin{eqnarray}
    \int_{\exp(\log x)^{1-\frac{\delta}{2}}}^{\infty} \frac{d \Psi(t,p)}{t} \ll \exp(-(\log x)^{\frac{\delta}{2}}) .
\end{eqnarray}
Hence, combining (5.13) with (5.16), we get
\begin{eqnarray}
    \mathcal{S}_{2} \ll \frac{x(\log \log x)^{k-2} }{\exp((\log x)^{\frac{\delta}{2}})p\log p}.
\end{eqnarray}
Thus, using (5.11), (5.12), and (5.17), we get that
\begin{eqnarray}
    \Psi_k^{(p)}(x) \ll \mathcal{S}_{1} + \mathcal{S}_{2} \ll \frac{x (\log \log x)^{k-2} \log p}{p\log x} + \frac{x(\log \log x)^{k-2} }{\exp((\log x)^{\frac{\delta}{2}})p\log p}. \nonumber
\end{eqnarray}
Next, by definition, we have that
\begin{align}
    \Psi_k(x,T) = \sum_{p \leq T} \Psi_k^{(p)}(x) &\ll \frac{x(\log \log x)^{k-2} }{\log x} \sum_{p \leq T} \frac{\log p}{p} +\frac{x(\log \log x)^{k-2} }{\exp(\log x)^{\frac{\delta}{2}}} \sum_{p \leq T} \frac{1}{p \log p} \nonumber  \\
    &\ll \frac{x(\log \log x)^{k-2} }{\log x}\log T + \frac{x(\log \log x)^{k-2} }{\exp(\log x)^{\frac{\delta}{2}}}\nonumber \\
    &\ll \frac{x(\log \log x)^{k-2}  \log T}{\log x} .
\end{align}
Hence, we have our required result.
 \qed \\ \\
Now, following the same trend as above, we move on to prove the uniform distribution of the $k^{th}$ largest prime factor of $n$ in reduced residue classes $mod\;\ell$. \\ \\
\textbf{Theorem 5.3:} \textit{For each integer positive integer $k$ and $\ell \geq 2$, the sequence $\{P_k(n)\}_n$ of the $k^{th}$ largest prime factors of $n$, is uniformly distributed in the reduced residue classes modulo $\ell$. More precisely, for each fixed $\ell \geq 2$, and any $1 \leq j \leq \ell$ satisfying $(j,\ell)=1$, we have that}
\begin{eqnarray}
    N_k(x;\ell,j):=\sum_{\substack{n\leq x\\ P_k(n) \equiv j\;(mod\;l)}} 1 = \frac{x}{\varphi(\ell)}+O\left(\frac{x(\log\log x)^k}{\log x}\right).
\end{eqnarray}
\textbf{Proof:} We start by looking at integers $n$, such that $\omega(n)\leq k-1$. So for all such $n \leq x$, we get due to Landau that
\begin{eqnarray}
    \sum_{\substack{n \leq x\\ \omega(n)\leq k-1}} 1 = \frac{x(\log\log x)^{k-2}}{(k-2)! \log x} + O_k\left(\frac{x(\log\log x)^{k-3}}{(k-3)! \log x}\right) .
\end{eqnarray}
For a fixed prime $p$, let $S_k(x,p)$ denote the set of integers $\leq x$ with at least $k$ distinct prime factors where the $i^{th}$ largest prime factors do not repeat, for $1 \leq i \leq k-1$, and $p$ is the $k^{th}$ largest prime factor. Using (5.20) and \textit{Lemma 5.1}, we have that
\begin{align}
    \sum_{p \leq x^{\frac{1}{k}}} |S_k(x,p)|&=[x]-\sum_{\substack{n \leq x\\ \omega(n)\leq k-1}} 1 + O_k\left(\frac{x(\log\log x)^{k-2}}{\log x}\right) \nonumber \\&= x+ O_k\left(\frac{x(\log\log x)^{k-2}}{\log x}\right).
\end{align}
We choose $N \in S_k(x,p)$, so $N=m.pp_{k-1}p_{k-2}\cdots  p_1$, where $P_i(N)=p_i$ and $P_i(N)$ is the $i^{th}$ largest prime factor of $N$, for $1 \leq i \leq k-1$. Thus $P_k(N)=p$, $P_1(m) \leq p$, and $p<p_{k-1}<\cdots  <p_1$. Then, with $q_j$ denoting the $j$-th prime number, we have 
\begin{align}
    |S_k(x,p)| &= \sum_{\substack{m<\frac{x}{p^k} \\ P_1(m)\leq p}}\;\;\sum_{p<p_{k-1}<\left(\frac{x}{mp}\right)^{\frac{1}{k-1}}}\;\;\sum_{p_{k-1}<p_{k-2}<\left(\frac{x}{mpp_{k-1}}\right)^{\frac{1}{k-2}}}\;\cdots  \;\sum_{p_2<p_1<\frac{x}{mpp_{k-1}\cdots   p_2}}1 \nonumber \\
    &=  \sum_{p<p_1<\frac{x}{2\cdot3\cdots   q_{k-1}}}\;\;\sum_{p<p_2<p_1}\;\;\sum_{p<p_3<p_2}\;\;\cdots  \;\;\sum_{p<p_{k-1}<p_{k-2}}\;\;\sum_{\substack{m \leq \frac{x}{pp_{k-1}\cdots  p_1}\\ P_1(m)\leq p}} 1 \nonumber \\
    &=  \sum_{p<p_1<\frac{x}{2\cdot3\cdots   q_{k-1}}}\;\;\sum_{p<p_2<p_1}\;\;\sum_{p<p_3<p_2}\;\;\cdots  \;\;\sum_{p<p_{k-1}<p_{k-2}} \Psi\left(\frac{x}{pp_{k-1}\cdots  p_1};p\right) .
\end{align}
Therefore, from (5.21) and (5.22), we obtain that
\begin{align}
    \sum_{p \leq x^{\frac{1}{k}}}|S_k(x,p)| &= \sum_{p\leq x^{\frac{1}{k}}}\;\;\sum_{p<p_1<\frac{x}{2\cdot3\cdots   q_{k-1}}}\;\;\sum_{p<p_2<p_1}\;\;\cdots   \;\; \sum_{p<p_{k-1}<p_{k-2}}\Psi\left(\frac{x}{pp_{k-1}\cdots   p_1};p\right) \nonumber \\
    &=x+O_k\left(\frac{x(\log\log x)^{k-2}}{\log x}\right).
\end{align}
Analogous to (3.18), we need to rewrite the multiple sum on the left in (5.23)
with inner sum over $p$, and this is
\begin{align}
    \sum_{p \leq x^{\frac{1}{k}}}|S_k(x,p)| &= \sum_{p_1<\frac{x}{2\cdot3\cdots   q_{k-1}}}\sum^*_{p_2}\sum^*_{p_3}\cdots  \sum^*_p \Psi\left(\frac{x}{pp_{k-1}\cdots   p_1};p\right) \nonumber \\
    &=x+O_k\left(\frac{x(\log\log x)^{k-2}}{\log x}\right),
\end{align}
where the $*$ over the summations mean that the conditions on the primes are
\begin{align}
    p\leq x^{\frac{1}{k}};\quad p<p_{k-1}<p_{k-2}<\cdots  <p_1\leq \frac{x}{2\cdot3\cdot\cdots  \cdot q_{k-1}} ;\quad\text{and}\quad pp_{k-1}\cdots  p_2p_1\leq  x. 
\end{align}
We note here that \textit{Theorem 5.2} gives us
\begin{eqnarray}
    \sum_{p \leq T} |S_k(x,p)| = O\left(\frac{x(\log\log x)^{k-2}\log T}{\log x} \right), 
\end{eqnarray}
where $T \leq e^{(\log x)^{1-\delta}}$, for some $\delta>0$.
So, we simply need to consider the summation $\sum_{T < p \leq x^{\frac{1}{k}}} |S_k(x;p)|$ for our required estimates. We have from (5.24) and (5.26) that
\begin{align}
    \sum_{T<p\leq x^{\frac{1}{k}}}|S_k(x,p)| &= \sum_{T<p_1\leq \frac{x}{T^{k-1}}}\sum_{T<p_2<p_1}^*\sum_{T<p_3<p_2}^*\cdots  \sum_{T<p<p_{k-1}}^*\Psi\left(\frac{x}{pp_1\cdots   p_{k-1}};p\right) \nonumber \\
    &= x+O_k\left(\frac{x(\log\log x)^{k-2}\log T}{\log x}\right),
\end{align} 
where the * on the top of each summation represents the conditions in (5.25), but we now also have $p>T$.
%\begin{align}
 %   T<p\leq x^{\frac{1}{k}};\quad p<p_{k-1}<p_{k-2}<\cdots  <p_1\leq \frac{x}{2\cdot3\cdot\cdots  \cdot P_{k-1}} ;\quad\text{and}\quad pp_{k-1}\cdots  p_2p_1\leq  x. \nonumber
%\end{align}
To evaluate this multi-sum, we split it into $k$ terms $\Sigma^T_j$, $j=1,2,\cdots  , k$, extending what we did in (3.21)-(3.24), and compare each sum with the corresponding integrals $I_j^T$, and evaluate the resulting error. That is, we write
%each with $k$-summations and denote each of them with $\Sigma_i^T$, $1 \leq i \leq k$. We write from (5.25) that
\begin{align}
     \sum_{T<p\leq x^{\frac{1}{k}}}|S_k(x,p)| = \Sigma_1^T+\Sigma_2^T+\cdots  +\Sigma_k^T  = x+O\left(\frac{x(\log\log x)^{k-2}\log T}{\log x}\right),
\end{align}
where the sums $\Sigma^T_j$ are generated using the following refined
conditions: If $p_i<x^{1/k}$, then for each $p_j$ with $j>k$, we automatically have $p_j<x^{1/k}$ because $p_j<p_i$. So we consider two cases, namely $p_j<x^{1/k}$ and $x^{1/k}<p_j<p_{j-1}$ in succession to generate the $k$ sums $\Sigma^T_i$ as
follows:
\begin{align}
    \Sigma_1^T&:= \sum_{T<p_1\leq x^{\frac{1}{k}}}\;\;\sum_{T<p_2<p_1}\;\;\cdots  \;\;\sum_{T<p<p_{k-1}} \Psi\left(\frac{x}{pp_1\cdots   p_{k-1}};p\right) \nonumber \\
        \Sigma_i^T&:= \sum_{x^{\frac{1}{k}}<p_1\leq \frac{x}{2\cdot3\cdots   q_{k-1}}}\
        \cdots   \sum_{x^{\frac{1}{k}}<p_{i-1}\leq p_i}\;\;\sum_{T<p_i\leq x^{\frac{1}{k}} }\;\;\sum_{T<p_{i+1}<p_i}\cdots  \sum_{T<p<\min\left(p_{k-1},\frac{x}{p_{k-1}\cdots   p_1}\right)} \Psi\left(\frac{x}{pp_{k-1}\cdots   p_1};p\right),\nonumber \\ &\hspace{12.5cm}\text{when $2 \leq i \leq k-1$};
        \nonumber \\
    \Sigma_{k}^T&:= \sum_{x^{\frac{1}{k}}<p_1\leq \frac{x}{2\cdot3\cdots   q_{k-1}}}\;\;\cdots   \;\;\sum_{x^{\frac{1}{k}}< p_{k-1}\leq p_{k-2}}\;\;\sum_{T<p<\frac{x}{p_1p_2\cdots   p_{k-1}}} \Psi\left(\frac{x}{pp_{k-1}\cdots   p_1};p\right) . \nonumber 
\end{align}
As we already know from the proof of \textit{Theorem 3.4}, the key step now is to replace the innermost summations on the $k^{th}$ largest prime factor $p$ in each of the multi-sums $\Sigma_i^T$ using corresponding integrals. We define pseudo-integrals $\mathcal{I}_i^T$, for $1 \leq i \leq k$ as follows:
\begin{align}
  \mathcal{I}_1^T&:= \sum_{T<p_1\leq x^{\frac{1}{k}}}\;\;\sum_{T<p_2<p_1}\;\;\cdots  \;\;\int_{T}^{p_{k-1}} \Psi\left(\frac{x}{tp_{k-1}\cdots  p_1};t\right)\frac{dt}{\log t};  \nonumber \\
    \mathcal{I}_i^T&:= \sum_{x^{\frac{1}{k}}<p_1\leq \frac{x}{2\cdot3\cdots   q_{k-1}}}\cdots  \sum_{x^{\frac{1}{k}}<p_{i-1}\leq p_i}\;\;\sum_{T<p_i\leq x^{\frac{1}{k}} }\;\;\sum_{T<p_{i+1}<p_i}\cdots  \;\;\int_{T}^{\min\left(p_{k-1},\frac{x}{p_{k-1}\cdots  p_1}\right)} \Psi\left(\frac{x}{tp_{k-1}\cdots  p_1};t\right)\frac{dt}{\log t} ,\nonumber \\
    &\hspace{12.5cm}\text{when $2 \leq i \leq k-1$} ;\nonumber \\
    \mathcal{I}_{k}^T&:= \sum_{x^{\frac{1}{k}}<p_1\leq \frac{x}{2\cdot3\cdot\cdots  \cdot q_{k-1}}}\;\;\cdots  \;\;\sum_{x^{\frac{1}{k}}< p_{k-1}\leq p_{k-2}}\;\;\int_{T}^{\frac{x}{p_1p_2\cdots  p_{k-1}}} \Psi\left(\frac{x}{tp_{k-1}\cdots  p_1};t\right)\frac{dt}{\log t} .
\end{align}
This replacement gives us from (5.28) that
\begin{align}
    \mathcal{I}_1^T+\cdots  +\mathcal{I}_k^T+E_1^T+\cdots  +E_k^T = x+ O\left(\frac{x(\log\log x)^{k-2}\log T}{\log x}\right),
\end{align}
where the error terms $E_i=\Sigma^T_i-\mathcal{I}_i^T$, for $1\leq i\leq k$ can be written as
\begin{align}
  E_1 &=  \sum_{T<p_1\leq x^{\frac{1}{k}}}\;\;\sum_{T<p_2<p_1}\;\;\cdots  \;\;\sum_{n \leq \frac{x}{Tp_{k-1}\cdots  p_1}}\left(\sum_{\max(T,P_1(n))\leq p \leq \min\left(\frac{x}{np_{k-1}\cdots  p_1},p_{k-1}\right)} 1 - \int_{\max(T,P_1(n))}^{\min\left(\frac{x}{np_{k-1}\cdots  p_1},p_{k-1}\right)}\frac{dt}{\log t}            \right); \nonumber \\
    E_i &= \sum_{x^{\frac{1}{k}}<p_1\leq \frac{x}{2\cdot3\cdots q_{k-1}}}\cdots  \sum_{x^{\frac{1}{k}}<p_{i-1}\leq p_i}\;\;\sum_{T<p_i\leq x^{\frac{1}{k}} }\;\;\sum_{T<p_{i+1}<p_i}\;\;\cdots   \nonumber \\
    &\hspace{1cm} \sum_{n \leq \frac{x}{Tp_{k-1}\cdots  p_1}}\left(\sum_{\max(T,P_1(n))\leq p \leq \min\left(\frac{x}{np_{k-1}\cdots  p_1},p_{k-1}\right)} 1 - \int_{\max(T,P_1(n))}^{\min\left(\frac{x}{np_{k-1}\cdots  p_1},p_{k-1}\right)}\frac{dt}{\log t} \right),\;\text{when $2 \leq i \leq k-1$} ;\nonumber \\
    E_k&= \sum_{x^{\frac{1}{k}}<p_1\leq \frac{x}{2\cdot3\cdots q_{k-1}}}\;\;\cdots  \;\;\sum_{x^{\frac{1}{k}}< p_{k-1}\leq p_{k-2}}\;\;\sum_{n \leq \frac{x}{Tp_{k-1}\cdots  p_1}}\left(\sum_{\max(T,P_1(n))\leq p \leq \frac{x}{np_{k-1}\cdots  p_1}} 1 - \int_{\max(T,P_1(n))}^{\frac{x}{np_{k-1}\cdots  p_1}}\frac{dt}{\log t}            \right) .
\end{align}
In obtaining the above expressions for $E_i$, we have, as in (3.26)
\begin{align}
    \text{replaced} \quad \Psi\left(\frac{x}{pp_{k-1}p_{k-2}\cdots   p_1};p\right) \quad \text{with}\quad
    \sum_{\substack{n\le \frac{x}{pp_{k-1}p_{k-2}\cdots   p_1} \\ P_1(n)<p}}1   \nonumber
\end{align}
in the innermost sum, and
\begin{align}
  \text{replaced} \quad \Psi\left(\frac{x}{tp_{k-1}p_{k-2}\cdots   p_1};t\right) \quad \text{with}\quad
    \sum_{\substack{n\le \frac{x}{tp_{k-1}p_{k-2}\cdots   p_1} \\ P_1(n)<t}}1     \nonumber
\end{align}
in the integral, and interchanged the order of summation over $p$ and $n$, and interchanged the order of the integral over $t$ and the summation over $n$. So we are able to use the strong form of the Prime Number Theorem to get a bound on the innermost expression in (5.31) to get a bound on each of the $E_i$. Once again, in getting this bound, we use the fact that the error term $\frac{t}{e^{\sqrt{\log\,t}}}$ is an increasing function of $t$. So we just bound the error by using only $t=\frac{x}{p_1p_2\cdots   p_{k-1}}$ and the trivial bound $p_i\leq x$ to get
\begin{eqnarray}
    |E_i| \ll\sum_{p_1<x}\sum_{p_2<x}\cdots  \sum_{p_{k-1}<x}\sum_{n<x}\;\;\frac{x}{np_{k-1}p_{k-2}\cdots   p_1e^{\sqrt{\frac{x}{np_1p_2\cdots   p_{k-1}}}}}\ll \frac{x\log x(\log\log x)^{k-1}}{e^{\sqrt{\log T}}},
\end{eqnarray}
because
\begin{align}
\frac{x}{np_1p_2\cdots   p_{k-1}}\ge T.
\end{align}
Therefore, using (5.32) in (5.30), we get that
\begin{eqnarray}
    \mathcal{I}_1^T+\cdots  +\mathcal{I}_k^T = x + O\left(\frac{x\log x(\log\log x)^{k-1}}{\exp\{\sqrt{\log T}\}} \right)+O\left( \frac{x(\log\log x)^{k-2}\log T}{\log x} \right).
\end{eqnarray}
Hence, choosing $T = \exp\{(2\log\log x)^2\}$, we get that
\begin{eqnarray}
    \mathcal{I}_1^T+\cdots  +\mathcal{I}_k^T = x + O\left(\frac{x(\log\log x)^k}{\log x}\right) .
\end{eqnarray}
Now, fixing an arbitrary integer $\ell \geq 2$, for any integer $1 \leq j \leq \ell$ such that $(j,\ell)=1$, we consider the set $S_{k,\ell,j}(x)$ of integers $n \leq x$ with square-free $i^{th}$ largest primes, where $1 \leq i \leq k-1$, and further satisfying $\omega(n) \geq k$ and $P_k(n) \equiv j\;(mod\;\ell)$. Therefore, using \textit{Lemma 5.1} and a similar argument used to get (5.21), we have
\begin{eqnarray}
    N_k(x;\ell,j) = |S_{k,\ell,j}(x)| + O_k\left(\frac{x(\log\log x)^{k-2}}{\log x}\right).
\end{eqnarray}
Also, we note here that
\begin{eqnarray}
    |S_{k,\ell,j}(x)| = \sum_{\substack{p \leq x^{\frac{1}{k}}\\ p \equiv j\;(mod\;\ell)}} |S_k(x,p)|.
\end{eqnarray}
Of course, we know from \textit{Theorem 5.2} that
\begin{eqnarray}
    \sum_{\substack{p \leq T\\ p \equiv j\;(mod\;\ell)}} |S_k(x,p)|=O\left(\frac{x(\log\log x)^{k-2}\log T}{\log x}\right),
\end{eqnarray}
where $T\leq e^{(\log x)^{1-\delta}}$, for some $\delta>0$. Focusing on the sum $\sum_{T<p \leq x^{\frac{1}{k}};\;p \equiv j\;(mod\;\ell)} |S_k(x,p)|$, we split it into $k$ many terms, each with $k$ sums as done previously in (5.26). Note that each of the $k$ sums will look exactly similar to the ones we already had earlier, along with the residue class condition on $p$. We then have, from (5.36), (5.37), and (5.38) that 
\begin{align}
    N_k(x;\ell,j) = \Sigma_{1,\ell,j}^T+\cdots  +\Sigma_{k,\ell,j}^T + O\left(\frac{x(\log\log x)^{k-2}\log T}{\log x}\right) ,
\end{align}
where 
the terms $\Sigma_{i,\ell,j}^T$, for each $1 \leq i \leq k$, are defined as follows:
\begin{align}
  \Sigma_{1,\ell,j}^T&:= \sum_{T<p_1\leq x^{\frac{1}{k}}}\;\;\sum_{T<p_2<p_1}\;\;\cdots  \;\;\sum_{\substack{T<p<p_{k-1} \\ p \equiv j\;(mod\;\ell)}} \Psi\left(\frac{x}{pp_{k-1}\cdots  p_1};p\right); \nonumber \\
  \Sigma_{i,\ell,j}^T&:= \sum_{x^{\frac{1}{k}}<p_1\leq \frac{x}{2\cdot3\cdots  q_{k-1}}}\cdots  \sum_{x^{\frac{1}{k}}<p_{i-1}\leq p_i}\;\;\sum_{T<p_i\leq x^{\frac{1}{k}} }\;\;\sum_{T<p_{i+1}<p_i}\cdots  \sum_{\substack{T<p<\min(p_{k-1},\frac{x}{p_{k-1}\cdots  p_1})\\ p \equiv j\;(mod\;\ell)}} \Psi\left(\frac{x}{pp_{k-1}\cdots  p_1};p\right),\nonumber \\ &\hspace{12.5cm}\text{when $2 \leq i \leq k-1$} ;\nonumber \\
  &=\nonumber \\
    \Sigma_{k,\ell,j}^T&:= \sum_{x^{\frac{1}{k}}<p_1\leq \frac{x}{2\cdot3\cdots q_{k-1}}}\;\;\cdots  \;\;\sum_{x^{\frac{1}{k}}< p_{k-1}\leq p_{k-2}}\;\;\sum_{\substack{T<p<\frac{x}{p_1p_2\cdots  p_{k-1}} \\ p \equiv j\;(mod\;\ell)}} \Psi\left(\frac{x}{pp_{k-1}\cdots  p_1};p\right) . \nonumber
\end{align}
Similar to what we did earlier, we replace each of these multi-sums with the following expressions denoted by $\mathcal{I}_{i,\ell,j}^T$:
\begin{align}
  \mathcal{I}_{1,\ell,j}^T&:= \sum_{T<p_1\leq x^{\frac{1}{k}}}\;\;\sum_{T<p_2<p_1}\;\;\cdots  \;\;\frac{1}{\varphi(\ell)}\int_{T}^{p_{k-1}} \Psi\left(\frac{x}{tp_{k-1}\cdots  p_1};t\right)\frac{dt}{\log t} ; \nonumber \\
    \mathcal{I}_{i,\ell,j}^T&:= \sum_{x^{\frac{1}{k}}<p_1\leq \frac{x}{2\cdot3\cdots q_{k-1}}}\cdots  \sum_{x^{\frac{1}{k}}<p_{i-1}\leq p_i}\;\;\sum_{T<p_i\leq x^{\frac{1}{k}} }\;\;\sum_{T<p_{i+1}<p_i}\cdots  \;\;\frac{1}{\varphi(\ell)}\int_{T}^{\min(p_{k-1},\frac{x}{p_{k-1}\cdots  p_1})} \Psi\left(\frac{x}{tp_{k-1}\cdots  p_1};t\right)\frac{dt}{\log t} ,\nonumber \\
    &\hspace{12.5cm}\text{when $2 \leq i \leq k-1$}; \nonumber \\
    &=\nonumber \\
    \mathcal{I}_{k,\ell,j}^T&:= \sum_{x^{\frac{1}{k}}<p_1\leq \frac{x}{2\cdot3\cdots q_{k-1}}}\;\;\cdots  \;\;\sum_{x^{\frac{1}{k}}< p_{k-1}\leq p_{k-2}}\;\;\frac{1}{\varphi(\ell)}\int_{T}^{\frac{x}{p_1p_2\cdots  p_{k-1}}} \Psi\left(\frac{x}{tp_{k-1}\cdots  p_1};t\right)\frac{dt}{\log t} .\nonumber
\end{align}
Next we look at the differences $E_{i,\ell,j}:=\Sigma_{i,\ell,j}^T-\mathcal{I}_{i,\ell,j}^T$. Using the above replacements in (5.39), we have
\begin{align}
     N_k(x;\ell,j) = \mathcal{I}_{1,\ell,j}^T+\cdots  +\mathcal{I}_{k,\ell,j}^T + E_{1,\ell,j}+\cdots  +E_{k,\ell,j} +O\left(\frac{x(\log\log x)^{k-2}\log T}{\log x}\right),
\end{align}
where the error terms can be written as
\begin{align}
  E_{1,\ell,j} &=  \sum_{T<p_1\leq x^{\frac{1}{k}}}\;\;\sum_{T<p_2<p_1}\;\;\cdots  \;\;\sum_{n \leq \frac{x}{Tp_{k-1}\cdots  p_1}}\left(\sum_{\substack{\max(T,P_1(n))\leq p \leq \min(\frac{x}{np_{k-1}\cdots  p_1},p_{k-1})\\ p \equiv j\;(mod\;\ell)}} 1 - \frac{1}{\varphi(\ell)}\int_{\max(T,P_1(n))}^{\min(\frac{x}{np_{k-1}\cdots  p_1},p_{k-1})}\frac{dt}{\log t}            \right); \nonumber \\
    E_{i,\ell,j} &= \sum_{x^{\frac{1}{k}}<p_1\leq \frac{x}{2\cdot3  \cdots q_{k-1}}}\cdots  \sum_{x^{\frac{1}{k}}<p_{i-1}\leq p_i}\;\;\sum_{T<p_i\leq x^{\frac{1}{k}} }\;\;\sum_{T<p_{i+1}<p_i}\;\;\cdots   \nonumber \\
    &\hspace{1cm} \sum_{n \leq \frac{x}{Tp_{k-1}\cdots  p_1}}\left(\sum_{\substack{\max(T,P_1(n))\leq p \leq \min(\frac{x}{np_{k-1}\cdots  p_1},p_{k-1})\\ p \equiv j\;(mod\;\ell)}} 1 - \frac{1}{\varphi(\ell)}\int_{\max(T,P_1(n))}^{\min(\frac{x}{np_{k-1}\cdots  p_1},p_{k-1})}\frac{dt}{\log t} \right),\;\text{when $2 \leq i \leq k-1$}; \nonumber \\
    E_{k,\ell,j}&= \sum_{x^{\frac{1}{k}}<p_1\leq \frac{x}{2\cdot3\cdots q_{k-1}}}\;\;\cdots  \;\;\sum_{x^{\frac{1}{k}}< p_{k-1}\leq p_{k-2}}\;\;\sum_{n \leq \frac{x}{Tp_{k-1}\cdots  p_1}}\left(\sum_{\substack{\max(T,P_1(n))\leq p \leq \frac{x}{np_{k-1}\cdots  p_1}\\ p \equiv j\;(mod\;\ell)}} 1 - \frac{1}{\varphi(\ell)}\int_{\max(T,P_1(n))}^{\frac{x}{np_{k-1}\cdots  p_1}}\frac{dt}{\log t}            \right). \nonumber
\end{align}
Using the strong form of the Prime Number Theorem in Arithmetic Progressions and the summation estimates, we then deduce that for every $1 \leq i \leq k$,
\begin{eqnarray}
    |E_{i,\ell,j}|\ll \frac{x\log x(\log\log x)^{k-1}}{\exp\{\sqrt{\log T}\}} .
\end{eqnarray}
Hence, from (5.40) and (5.41), we get that
\begin{eqnarray}
    N_k(x;\ell,j) = \mathcal{I}_{1,\ell,j}^T+\cdots  +\mathcal{I}_{k,\ell,j}^T + O\left(\frac{x\log x(\log\log x)^{k-1}}{\exp\{\sqrt{\log T}\}} \right)+O\left( \frac{x(\log\log x)^{k-2}\log T}{\log x} \right).
\end{eqnarray}
But then, we note that $\mathcal{I}_{i,\ell,j}^T = \frac{1}{\varphi(\ell)}\mathcal{I}_i^T$, for each $1 \leq i\leq k$. Therefore, using this in (5.42), and choosing $T = \exp\{(2\log\log x)^2\}$, we finally get from (5.35) that
\begin{eqnarray}
    N_k(x;\ell,j) = \frac{1}{\varphi(\ell)}(\mathcal{I}_1^T+\cdots  +\mathcal{I}_k^T) + O\left(\frac{x(\log\log x)^k}{\log x}\right) = \frac{x}{\varphi(\ell)} + O\left(\frac{x(\log\log x)^k}{\log x}\right). \nonumber
\end{eqnarray}
Therefore, we have our required result that the sequence of $k^{th}$ largest prime factors of integers is uniformly distributed in the reduced residue classes. \qed \\ \\
As remarked at the end of \S3, it is quite fascinating to see that the averages of any $k^{th}$ largest prime factors, where $k$ is arbitrarily chosen, is the same in reduced residue classes, and every time it is $\frac{1}{\varphi(\ell)}$, when we consider the reduced residue classes $mod\;\ell$. This is an extremely important property, due to which our generalizations (see next section) of our previous results are possible. With unequal averages for every sequence $k^{th}$ largest prime factors, for different integers $k$, providing a generalization would have been acutely hard. This is another reason why the use of induction is so powerful in the final theorems that are to come in the next section.

\section{The Main Theorems for all orders $k\geq 4$}
In this section, we extend the results of \S4 to all orders $k\geq 4$. We start with the following:
 \\ \\
\textbf{Theorem 6.1} \textit{For any fixed integer $k \geq 4$, and integers $j,\ell$ satisfying $1 \leq j \leq \ell$ and $(j,\ell)=1$, we have that}
\begin{eqnarray}
    M_{\omega^{k-1}}(x,j,\ell) := \sum_{\substack{n \leq x \\ p_1(n) \equiv j\;(mod\;\ell)}} \mu(n)\omega(n)^{k-1} \ll \frac{x(\log\log x)^{k+2}}{\log x} .
\end{eqnarray}
\textit{Equivalently,}
\begin{align}
    M_{\omega_{k-1}}(x,j,\ell) := \sum_{\substack{n \leq x \\ p_1(n) \equiv j\;(mod\;\ell)}} \mu(n)\omega_{k-1}(n) \ll \frac{x(\log\log x)^{k+2}}{\log x},
\end{align}
\textit{where $\omega_{k-1}(n) = {\omega(n)-1 \choose k-1}$}.\\ \\
\textbf{Proof:} We start by recalling the function that we defined at the start of the proof of \textit{Theorem 4.1}: $f$ is an arithmetic function on primes defined as follows:
\begin{eqnarray}
    f(n) = \begin{cases}
        1, \;\text{if $n$ is a prime and $n \equiv j\;(mod\;\ell)$} \\
        0,\;\text{otherwise.}
    \end{cases}
\end{eqnarray}
Let us also recall the following version of the Duality Lemma in (1.11), that we will need in this proof:
\begin{eqnarray}
    \sum_{1<d|n} \mu(d){\omega(d)-1 \choose k-1}f(p_1(d)) = (-1)^kf(P_k(d)).
\end{eqnarray}
To prove \textit{Theorem 6.1}, we will first establish (6.2) for any given $k$ and utilize the truth of (6.2) for $t \leq k-1$ to deduce (6.1). \\ \\
Using M\"obius Inversion on (6.4), we get
\begin{eqnarray}
    \sum_{n \leq x}\mu(n)(\omega(n)-1)\cdots  (\omega(n)-k+1)f(p_1(n)) = (-1)^k(k-1)!\sum_{n \leq x}\sum_{d|n}\mu\left(\frac{n}{d}\right) f(P_k(d)).
\end{eqnarray}
Without worrying about the constant in front, we estimate the double sum on the RHS of (6.5) using the hyperbola method. We rewrite the double sum as follows:
\begin{eqnarray}
\sum_{n \leq x}\sum_{d|n}\mu\left(\frac{n}{d}\right) f(P_k(d)) = \sum_{m \leq T}\mu(m)\sum_{d \leq \frac{x}{m}} f(P_k(d)) + \sum_{d \leq \frac{x}{T}} f(P_k(d))  \sum_{T<m  \leq \frac{x}{d}} \mu(m),
\end{eqnarray}
for some suitable $T$ to be chosen later. Using the estimate of the sum of $\mu(m)$, we get that the second double sum in (6.6) has the following estimate:
\begin{eqnarray}
     \sum_{d \leq \frac{x}{T}}f(P_k(d))\sum_{T<m\leq \frac{x}{d}}\mu(m) \ll \sum_{d \leq \frac{x}{T}}f(P_k(d)) \frac{x}{de^{\sqrt{\log (x/d)}}}\ll \frac{x\log x}{e^{\sqrt{\log T}}} .
\end{eqnarray}
Further, using \textit{Theorem 5.3}, we have that
\begin{eqnarray}
    \sum_{d \leq \frac{x}{m}} f(P_k(d)) = \frac{x}{m\varphi(\ell)} + O\left(\frac{x(\log\log x)^k}{m\log(x/m)}\right) .
\end{eqnarray}
Therefore, using (6.8) in the first double sum of (6.6), we have
\begin{align}
     \sum_{m \leq T}\mu(m)\sum_{d \leq \frac{x}{m}} f(P_t(d)) &= \frac{x}{\varphi(\ell)} \sum_{m \leq T}\frac{\mu(m)}{m}+O\left(x(\log\log x)^k\sum_{m\leq T}\frac{1}{m\log(x/m)}\right)\nonumber \\ 
     &\ll \frac{x}{\varphi(\ell)e^{\sqrt{\log (T)}}}+O\left(\frac{x\log T(\log\log x)^t}{\log (x/T)}\right) .
\end{align}
Therefore, choosing $T = e^{(2\log\log x)^2}$, we have from (6.6)-(6.9) that
\begin{eqnarray}
    \sum_{n \leq x}\sum_{d|n}\mu\left(\frac{n}{d}\right) f(P_k(d)) \ll \frac{x(\log\log x)^{k+2}}{\log x},
\end{eqnarray}
and hence, (6.2) follows from (6.5) and (6.10). Next, write
\begin{align}
    \omega(n)^{k-1} = \sum_{t=0}^{k-1} s_{k-1,t}\omega_{t}(n),
\end{align}
where $s_{k-1,t}$ are given in terms of Sterling numbers of the second kind, with $s_{k-1,k-1}=1$ (for the definition of Sterling numbers of the second kind, see 6.23 below). Hence, 
\begin{align}
    M_{\omega^{k-1}}(x;j,\ell) = \sum_{t=0}^{k-1}s_{k-1,t}M_{\omega_{t}}(x;j,\ell).
\end{align}
Then (6.1) follows by utilizing in (6.12) the upper bound for
$M_{\omega_t}(x;j,\ell)$  for all $t\le k-1$. That proves \textit{Theorem 6.1}.
%\begin{eqnarray}
 %   M_{\omega_{t-1}}(x,j,\ell) = \sum_{\substack{n \leq x \\p_1(n) \equiv j\;(mod \ell)}} \mu(n)(\omega(n)-1)\cdots  (\omega(n)-t+1)+(t-1) \sum_{\substack{n \leq x \\p_1(n) \equiv j\;(mod \ell)}} \mu(n)(\omega(n)-1)\cdots  (\omega(n)-t+2) \nonumber \\
  %  = \sum_{\substack{n \leq x \\p_1(n) \equiv j\;(mod \ell)}} \mu(n)(\omega(n)-1)\cdots  (\omega(n)-t+1) + (t-1) \sum_{\substack{n \leq x \\p_1(n) \equiv j\;(mod \ell)}} \mu(n)\omega(n)\cdots  (\omega(n)-t+3)\hspace{1cm} \nonumber \\
%    -(t-1)(t-2)\sum_{\substack{n \leq x \\p_1(n) \equiv j\;(mod \ell)}} \mu(n)(\omega(n)-1)\cdots  (\omega(n)-t+3) \hspace{1.1cm} \nonumber \\
 %   = \sum_{\substack{n \leq x \\p_1(n) \equiv j\;(mod \ell)}} \mu(n)(\omega(n)-1)\cdots  (\omega(n)-t+1) + (t-1)M_{\omega_{t-2}}(x,j,\ell)-(t-1)(t-2)M_{\omega_{t-3}}(x,j,\ell)\hspace{0.2cm}  \nonumber \\
  %  \cdots  +(-1)^{t}(t-1)!\sum_{\substack{n \leq x \\p_1(n) \equiv j\;(mod \ell)}} \mu(n)\hspace{1cm} \nonumber 
%\end{eqnarray}
\qed \\ \\
\textbf{Remark:} (6.1) for $k=1$ follows from \cite{KA1977} and for $k=2$ is established in \cite{AJ24}, although the estimates were much stronger in those cases.\\ \\
As done earlier, we now estimate the sum with the same summand as that of above in \textit{Theorem 6.1}, with an additional weight of the fractional part $\left\{\frac{x}{n}\right\}$. \\ \\
\textbf{Theorem 6.2:} \textit{For any fixed integer $k \geq 4$, and integers $j,\ell$ satisfying $1 \leq j \leq \ell$ and $(j,\ell)=1$, we have that}
\begin{eqnarray}
     \sum_{\substack{n \leq x\\ p_1(n)\equiv j\;(mod\;\ell)}}\mu(n)\omega(n)^{k-1}\left\{\frac{x}{n}\right\}\ll\frac{x(\log\log x)^{k+\frac{1}{2}}}{\sqrt{\log x}}.
\end{eqnarray}
\textbf{Proof:} Just like in \textit{Theorem 4.2}, we use \textit{Theorem A} from \S2 to prove \textit{Theorem 6.2}. Here, we choose the sequence $\{a_n\}_n$ to be $\mu(n)\omega(n)^{k-1}f(p_1(n))$, where $f$ is chosen to be the function defined in (6.3). Then, by \textit{Theorem 6.1}, we have that
\begin{align}
    \sum_{n\leq x}a_n \ll \frac{x(\log\log x)^{k+2}}{\log x}.
\end{align}
If we choose $\eta(x)=\frac{(\log\log x)^{k+2}}{\log x}$, it satisfies the conditions in (2.34a) and (2.34b). Further, we have that
\begin{eqnarray}
    \sum_{n\leq x}|a_n| = \sum_{n\leq x}\omega(n)^{k-1} \ll x(\log\log x)^{k-1}. 
\end{eqnarray}
Therefore, the correct choice of $\beta(x)$ in \textit{Theorem A} is then $\beta(x)=(\log\log x)^{k-1}$ as it satisfies the conditions (2.36). Hence, by (2.37), \textit{Theorem 6.2} follows.
 \qed \\ \\
As done in every case throughout the paper, we now estimate the exact same summation in the above theorem with the fractional part of $\frac{x}{n}$ replaced by its integral part. We have the following theorem: \\ \\ 
\textbf{Theorem 6.3:} \textit{For any fixed integer $k \geq 4$, and integers $j, \ell$ satisfying $1 \leq j \leq \ell$ and $(j,\ell)=1$, we have that}
 \begin{eqnarray}
     \sum_{\substack{n \leq x\\ p_1(n)\equiv j\;(mod\;\ell)}}\mu(n)\omega(n)^{k-1}\left[\frac{x}{n}\right]\ll\frac{x(\log\log x)^{k}}{\log x}.
 \end{eqnarray}
\textbf{Proof:} Let us start by looking at the following summation:
\begin{align}
    \sum_{\substack{n \leq x \\ p_1(n) \equiv j\;(mod\;\ell)}} \mu(n)\omega(n)\omega_{t-2}(n) \left[\frac{x}{n}\right],\nonumber 
\end{align}
where $\omega_{t-2}(n) = {\omega(n)-1 \choose t-2}$, where $t \geq 5$. Using the definition (6.3) of the function $f$, we have from the above sum that
\begin{align}
    \sum_{\substack{n \leq x \\ p_1(n) \equiv j\;(mod\;\ell)}} \mu(n)\omega(n)\omega_{t-2}(n)\left[\frac{x}{n}\right] &= \sum_{n \leq x } \mu(n)\omega(n)\omega_{t-2}(n)f(p_1(n))\left[\frac{x}{n}\right] \nonumber  \\
    &= \sum_{n \leq x }\sum_{d|n} \mu(d)\omega(d)\omega_{t-2}(d)f(p_1(d)).
\end{align}
Now, we notice that
\begin{align}
    \omega(n)\omega_{t-2}(n) &= \frac{1}{(t-2)!} \cdot \omega(n)(\omega(n)-1)(\omega(n)-2)\cdots  (\omega(n)-t+2) \nonumber \\
    &=\frac{1}{(t-2)!}\left[(\omega(n)-1)\cdots  (\omega(n)-t+1)+(t-1)(\omega(n)-1)\cdots  (\omega(n)-t+2)\right] \nonumber \\
    &= \frac{1}{(t-2)!}\left[(t-1)!\cdot {\omega(n)-1 \choose t-1}+(t-1)\cdot(t-2)!\cdot{\omega(n)-1 \choose (t-1)-1}\right] \nonumber \\
    &= (t-1)\left[{\omega(n)-1 \choose t-1}+{\omega(n)-1 \choose (t-1)-1}\right].
\end{align}
Now, using (6.18) in the summand of the double sum in the RHS of (6.17), we get by splitting the summation that
\begin{align}
     \sum_{\substack{n \leq x \\ p_1(n) \equiv j\;(mod\;\ell)}} \mu(n)\omega(n)\omega_{t-2}(n)\left[\frac{x}{n}\right] = (t-1)\left[\sum_{n\leq x}\sum_{d|n}\mu(d){\omega(d)-1\choose t-1}f(p_1(d))\right.\hspace{2cm}\nonumber \\
     \left.+\sum_{n\leq x}\sum_{d|n}\mu(d){\omega(d)-1\choose (t-1)-1}f(p_1(d))\right] .
\end{align}
Applying the higher order duality identity in (6.4) for $k=t$ and $k=t-1$ respectively in the above two double sums, we get that
\begin{align}
     \sum_{\substack{n \leq x \\ p_1(n) \equiv j\;(mod\;\ell)}} \mu(n)\omega(n)\omega_{t-2}(n)\left[\frac{x}{n}\right] &= (t-1)\left[\sum_{n \leq x}(-1)^{t}f(P_{t}(n))+\sum_{n\leq x}(-1)^{t-1}f(P_{t-1}(n))\right] \nonumber \\
     &= (t-1)(-1)^t\left[\sum_{n\leq x}f(P_t(n))-\sum_{n\leq x}f(P_{t-1}(n))\right].
\end{align}
By \textit{Theorem 5.3}, for both the sums inside the brackets in the RHS of (6.20), we have a cancellation of the main terms to deduce that
\begin{align}
     \sum_{\substack{n \leq x \\ p_1(n) \equiv j\;(mod\;\ell)}} \mu(n)\omega(n)\omega_{t-2}(n)\left[\frac{x}{n}\right]  \ll  \frac{x(\log\log x)^{t}}{\log x}.
\end{align}
It is important to note here that (6.21) has been proved in \cite{KA1977} and \cite{AJ24} for both $t=2,3$, where the error estimates are stronger than what has been proved above for $t=4$ (\textit{Theorem 4.3}) and $t \geq 4$. Now, we observe that
\begin{align}
    \omega(n)^{k-1} = \sum_{t=2}^{k} C_{t,k-1}\omega(n)\omega_{t-2}(n),
\end{align}
where $C_{t,k-1}$ are constants. Therefore, (6.21) and (6.22) together prove our required result.\qed \\ \\
\textbf{Remark:} Note that the case of $k=3$ corresponding to \textit{Theorem 6.3} is \textit{Theorem 4.3} in \S4. But we have used different methods to deduce these two results. In the proof of \textit{Theorem 4.3}, we have utilized a cancellation, namely $1-3+2=0$, which eliminates the main term of size $\frac{x}{\varphi(\ell)}$, thereby giving a bound $o(x)$. This cancellation is due to a property of
the Stirling numbers of the second kind which are defined by
\begin{align}
x^{k-1}=\sum^{k}_{j=1}S_{k,j}(x-1)(x-2)\cdots  (x-j+1).
\end{align}
It is known that the Stirling numbers of the second kind satisfy the following
recurrence:
\begin{align}
    S_{k,j} = jS_{k-1,j}+S_{k-1,j-1}, 
\end{align}
where $S_{k,j}$ is the member in the Stirling's pyramid in the $j^{th}$ position of the $k^{th}$ row. We set $S_{k,0}=0=S_{k,k+1}$ and note that $S_{k,1}=1=S_{k,k}$. We now claim that the following cancellation holds: \\ \\
\textbf{Lemma 6.4:}
\begin{align}
    \sum_{j=1}^k(-1)^{j-1}(j-1)!S_{k,j}=0.\nonumber
\end{align}
\textbf{Remark:} This is the cancellation property that generalizes what he had in (4.19) and
can be used to prove \textit{Theorem 6.3}. \\ \\
\textbf{Proof:} We observe that 
\begin{align}
     \sum_{j=1}^k(-1)^{j-1}(j-1)!S_{k,j} &= \sum_{j=1}^k(-1)^{j-1}(j-1)!\left\{jS_{k-1,j}+S_{k-1,j-1}\right\} \nonumber \\
     &= -\sum_{j=1}^{k}(-1)^jj!S_{k-1,j}+\sum_{j=0}^{k-1}(-1)^jj!S_{k-1,j}. \nonumber
\end{align}
Using the fact that $S_{k-1,k}=0=S_{k-1,0}$, we see that
\begin{align}
     \sum_{j=1}^k(-1)^{j-1}(j-1)!S_{k,j} = -\sum_{j=1}^{k-1}(-1)^jj!S_{k-1,j}+\sum_{j=1}^{k-1}(-1)^jj!S_{k-1,j}=0, \nonumber
\end{align}
which proves the lemma. \qed\\ \\
To prove \textit{Theorem 6.3} utilizing this cancellation, we need to write
\begin{align}
\omega^{k-1}(n)=\sum^{k}_{j=1}S_{k,j}(\omega(n)-1)(\omega(n)-2)\cdots  (\omega(n)-j+1) \nonumber
\end{align}
and utilize bounds for each of the sums $M_{\omega_t}(x;j,\ell)$ given by (6.2). The main term of size $\frac{x}{\varphi(\ell)}$ will be eliminated because of the cancellation given above, thereby leading to the bound for $M_{\omega^{k-1}}(x;j,\ell)$ given
in Theorem 6.3. The proof of \textit{Theorem 6.3} given above uses a different technique to bypass the use of the sum over $k$ terms by the introduction of the new term $\omega(n)\omega_{t-2}(n)$, which only requires an estimation of two double sums (see (6.18), (6.19) above).
%whereas an approach involving Stirling's formula above would give us $k+1$ many double sums.
\\ \\
We will now state our final theorem of this section, which, as seen earlier in every case, will be a culmination of \textit{Theorem 6.2} and \textit{Theorem 6.3}. \\ \\
\textbf{Theorem 6.5:} \textit{For any fixed integer $k \geq 4$, and integers satisfying $1 \leq j \leq \ell$ and $(j,\ell)=1$, we have that} 
\begin{eqnarray}
    \sum_{\substack{n \leq x \\ p_1(n)\equiv j \;(mod\;\ell)}} \frac{\mu(n)\omega(n)^{k-1}}{n} \ll \frac{(\log\log x)^{k+\frac{1}{2}}}{\sqrt{\log x}}.
\end{eqnarray} 
\textit{Further, taking $x \to \infty$ we have that}
\begin{eqnarray}
     \sum_{\substack{n=2 \\ p_1(n)\equiv j \;(mod\;\ell)}}^{\infty} \frac{\mu(n)\omega(n)^{k-1}}{n} = 0 .\nonumber
\end{eqnarray}
\textbf{Proof:} Combining the estimates (6.13) and (6.16), we have that
\begin{eqnarray}
    x  \sum_{\substack{n \leq x \\ p_1(n)\equiv j \;(mod\;\ell)}} \frac{\mu(n)\omega(n)^{k-1}}{n} \ll \frac{x(\log\log x)^{k+\frac{1}{2}}}{\sqrt{\log x}}.
\end{eqnarray}
Then, canceling $x$ on both sides of (6.26), we get (6.25) and hence, our main result is proved. \qed  \\ \\
We state the arithmetic density version of \textit{Theorem 6.5} as its corollary below. \\ \\
\textbf{Corollary 6.6:} \textit{For any fixed integer $k \geq 4$, and integers satisfying $1 \leq j \leq \ell$ and $(j,\ell)=1$, we have that} 
\begin{eqnarray}
     (-1)^k\sum_{\substack{n =1\\ p_1(n)\equiv j \;(mod\;\ell)}}^{\infty}  \frac{\mu(n){\omega(n)-1 \choose k-1}}{n} = \frac{1}{\varphi(\ell)}.
\end{eqnarray}

\section{Concluding Remarks}
In this paper, we have focused on the function $f(p)$ being the characteristic function of primes in an arithmetic progression. We conclude here by briefly discussing the more general situation when $f$ is a bounded function on the primes.
We also report the work of several authors on extensions of Alladi's duality to various algebraic settings.
\subsection{The Case of General $f$}
As noted in \S1, it was established by the first author in 1977 that if $f$ is a bounded function on the primes, and has an average (asymptotically) on the largest prime factor sequence $P_1(n)$ as in (1.4), then (1.5) holds and vice-versa. When $f$ is the characteristic function of primes in an arithmetic progression $j(mod\,\ell)$, then using the strong form of Prime Number Theorem for Arithmetic Progressions, Alladi \cite{KA1977} established (1.4) in quantitative form with $c=\frac{1}{\varphi(\ell)}$. Hence, (1.5) in quantitative form follows. \\ \\
As noted also in \S1, if $f$ is a bounded function on the primes such that asymptotically $f(P_1(n))$ and $f(P_2(n))$ have the same average $\kappa$, then (1.16) holds (see \cite{AJ24}). When $f$ is the characteristic functions of primes in an arithmetic progression $j(mod\,\ell)$, it was shown in \cite{AJ24}, that $f(P_(n))$
and $f(P_2(n))$ have the same average $\frac{1}{\varphi(\ell)}$, and so (1.18) holds. The deduction of (1.5) from (1.4) and of (1.16) from (1.14) and (1.15) makes use of the First Order Duality identity (1.3) and the Second Order Duality identity (1.10) when $k=2$. The more general situation is as follows: Let $k$ be a positive integer. If $f$ is a bounded function on the primes such that for $k\geq 2$,
\begin{align}
\lim_{x\to\infty}\frac{1}{x}\sum_{n\le x}f(P_j(n))=\kappa, \quad \text{for} \quad
j=1,2,\cdots  , k,    
\end{align}
then
\begin{align}
\sum^{\infty}_{n=2}\frac{\mu(n)\omega^{k-1}(n)f(p_1(n))}{n}=0.    
\end{align}
In order to deduce (7.2) from (7.1), we need to show that if $f$ is any bounded function on the primes, then for each positive integer $k$, 
\begin{align}
M_{\omega^{k-1}}(f, x):=\sum_{n\le x}\mu(n)\omega^{k-1}(n)f(p_1(n))=o(x).    
\end{align}
In \cite{KA1977}, (7.3) is established for $k=1$, and improved in \cite{AL82} to a best possible quantitative version. A quantitative version of (7.2) when
$k=2$ will be the subject of the paper by Alamoudi and Alladi \cite{Alladi-Alamoudi}. The general case
$k\ge 3$ of (7.3) will be treated in a subsequent paper of ours. 

\subsection{Algebraic Analogues of Duality between Prime Factors}
Forty years after the paper of the first author, spurred by a suggestion of Ken
Ono, Dawsey \cite{Da17} extended the results in \cite{KA1977} to a number field setting.
More precisely, consider a Galois extension $K$ of the rationals $\mathbb{Q}$, and $\mathcal{O}_K$,
the ring of integers in $K$. If $p$ is a rational integer prime, let $\mathfrak{p}$ be the
prime ideal contained in $\mathcal{O}_K$ which lies above $p$. If $p$ is unramified in $K$, let
$\left[\frac{K/\mathbb{Q}}{p}\right]$ denote the Artin symbol. Now consider any conjugacy class $C$ in the Galois group $G=Gal(G/Q)$. Then Dawsey's result is
\begin{align}
-\sum_{\substack{n\ge 2 \\ \left[\frac{K/\mathbb{Q}}{p_1(n)}\right]=C}}\frac{\mu(n)}{n}=\frac{|C|}{|G|}.    
\end{align}
When $K$ is an $\ell^{th}-$cyclotomic extension of $\mathbb{Q}$, $G$ can be identified with $\mathbb{Z}^*_{\ell}$, the set of reduced residues $(mod\,\ell)$ for some positive integer $\ell$.
Since $\mathbb{Z}^*_{\ell}$ is Abelian, each conjugacy class has just one element, and so
\begin{align}
\frac{|C|}{|G}=\frac{1}{\phi(\ell)}. \nonumber    
\end{align}
Just as Alladi used the strong form of the Prime Number Theorem for Arithmetic
Progressions, Dawsey used the strong form of the Chebotarev Density Theorem due to Lagarias and Odlyzko \cite{LO77} to show
\begin{align}
\sum_{\substack{2\le n\le x \\ \left[\frac{K/\mathbb{Q}}{P_1(n)}\right]=C}}1\sim \frac{|C|}{|G|}x,    
\end{align}
Then by methods involving Duality as in \cite{KA1977}, Dawsey deduces (7.4) from
(7.5).   \\ \\
Dawsey's paper led to a flurry of activity on algebraic extensions of the Duality results in \cite{KA1977}. Sweeting and Woo \cite{SW19} extending Dawsey's results to the case of $L$ being a Galois extension of a general number field $K$, and in that process generalized the M\"obius function by defining it in terms of products of prime ideals instead of products of primes. The results of Sweeting-Woo were further generalized by Kural et-al \cite{KMS}. There was also an extension of the results of \cite{KA1977}
by Wang \cite{WRam}, who replaced the Möbius function by the Ramanujan sum $c_m(n)$, which is the sum of the $m$-th powers of the primitive $n$-th roots of unity. Note that $c_1(n)=\mu(n)$. Subsequently, Duan, Wang, and Yi \cite{DWY}, discussed Alladi's Duality in global function fields. There are also other works of Wang on the first order duality - \cite{Wana}, in which $\mu(n)$ is replaced by Liouville's function $\lambda(n)$, and \cite{Wlog}, where $\mu(n)$ is replaced by $\mu(n)\log n$. \\ \\
While all these algebraic extensions of Alladi's First Order Duality were being
established, consequences of the Higher Order Duality identities in \cite{KA1977} remained unexplored until Alladi-Johnson in 2022 explored the Second Order Duality identity (1.10) for $k=2$ (see \cite{AJ24}), and established among other things, their
main result (1.12) in quantitative form. Motivated by this and Dawsey's theorem 
(7.4), Sengupta \cite{Se25} extended (1.12) to an algebraic setting, namely in the case of a Galois extension $K/\mathbb{Q}$. More recently, Jagannath Sahoo and Prasanna Nand Jha \cite{NJS} have extended
Alladi's Higher Order Duality Identities (for all $k$) to the function field setting, and
proved analogues of Alladi-Johnson's result (1.12) in the function field setting, which is the case of Second Order Duality ($k=2$). Similarly, Sengupta \cite{Se26}
has most recently  extended the results of Sweeting and Woo by establishing Higher Order Duality Identities (for all $k$) to the situation of $L$ being an
Galois extension of a number field $K$, and obtained analogues of (1.12),
namely the case $k=2$ in the case of Galois extensions $L/K$.  In view of the results of this paper, the following problems naturally occur for future exploration: \textit{(i)} to extend the results of sections 4 and 6 in the case of Galois extensions $L$ of a number field $K$ (a paper on this by Sengupta is in preparation), and \newline\textit{(ii)} to extend the results of sections 4 and 6 to the function field setting.

\bibliographystyle{alpha} % Styles: plain, alpha, IEEEtran, etc.
\bibliography{ref} % Points to references.bib

\vspace{1cm}

\noindent Krishnaswami Alladi: Department of Mathematics, University of Florida, Gainesville, FL 32611, United States \\
\textit{Email address:} alladik@ufl.edu\\ \\
Sroyon Sengupta: Department of Mathematics, University of Florida, Gainesville, FL 32611, United States \\
\textit{Email address:} sroyonsengupta1947@gmail.com

\end{document}